\documentclass[matzei,numbook]{svjour}
\usepackage{amsmath,amsfonts,amssymb}

\spnewtheorem*{Theorem}{Theorem}{\bf}{\it}
\spnewtheorem*{Proposition}{Proposition}{\bf}{\it}
\spnewtheorem*{Corollary}{Corollary}{\bf}{\it}
\spnewtheorem*{Lemma}{Lemma}{\bf}{\it}

\def\AMN{{\rm A}_{NM}}

\def\Bla{{\rm B}_{\la\lat}}
\def\Bll{{\rm B}_{l\ts\lt}}
\def\Bmm{{\rm B}_{\ts m\mt}}
\def\Bnn{{\rm B}_{\ts n\nt}}
\def\Bt{\ts\widetilde{\ns B\ts}}

\def\bx{{\boxed{\phantom{\square}}\kern-.4pt}}

\def\CC{{\mathbb C}}
\def\com{\ts,\hskip-.5pt}

\def\de{\delta}
\def\De{\Delta}

\def\End{\operatorname{End}\ts}

\def\FLa{F_{\ns\La\Lat}}
\def\FOm{F_{\Om\Omt}}

\def\ga{\gamma}
\def\gap{\gamma^{\,\prime}}
\def\ge{\geqslant}
\def\glL{{\mathfrak{gl}_L}}
\def\glM{{\mathfrak{gl}_M}}
\def\glMN{{\mathfrak{gl}_{\ts N+M}}}
\def\glN{{\mathfrak{gl}_N}}
\def\GOm{G_{\Om\Omt}}
\def\GOmp{G_{\Om\Omt}^{\,\,\prime}}

\def\id{{\rm id}}
\def\Ill{{\rm I}_{\ts l\ts\lt}}
\def\Imm{{\rm I}_{\ts m\mt}}

\def\io{\iota}

\def\kt{{\widetilde{k}}}

\def\la{\lambda}
\def\La{\Lambda}
\def\lap{\la^{\ts\prime}}
\def\Lap{\La^{\ts\prime}}
\def\Lapp{\La^{\ts\prime\prime}}
\def\lat{\widetilde{\la}}
\def\Lat{\tilde{\La}}
\def\Latp{\Lat^{\ts\prime}}
\def\Latpp{\Lat^{\ts\prime\prime}}
\def\lcd{\ts,\,\ldots,}
\def\le{\leqslant}
\def\lm{{\la/\ns\mu}}
\def\lmt{{\lat/\mut}}
\def\lt{{\tilde{l}}}

\def\mi{{\raise.5pt\hbox{-}}}
\def\mt{{\widetilde{m}}}
\def\mup{\mu^{\ts\prime}}
\def\mut{\widetilde{\mu}}
\def\mw{\kern-81pt}

\def\ns{\hskip-1pt}
\def\nt{{\widetilde{n}}}
\def\nup{\nu^{\ts\prime}}
\def\nut{\widetilde{\nu}}

\def\om{\omega}
\def\Om{\Omega}
\def\Omt{{\widetilde{\Om}}}
\def\op{\oplus}
\def\ot{\otimes}

\def\ph{\varphi}
\def\Pp{P^{\vee}}
\def\Pt{\ts\widetilde{\ns P\ts}}

\def\Qp{Q^{\,\vee}}

\def\Ra{R^{\hskip1pt\ast}}
\def\Rb{\,\overline{\hskip-2.5pt R\hskip-4.5pt\phantom{\bar{t}}}\ts}
\def\Rbp{\Rb^{\,\vee}}
\def\Rp{R^{\,\vee}}
\def\RR{{\mathbb R}}
\def\Rt{\ts\widetilde{\ns R\ts}}

\def\si{\sigma}

\def\Sll{S_{\ts\lt\ts l}}

\def\th{\theta}
\def\Th{\Theta}
\def\ts{\hskip1pt}

\def\UL{\operatorname{U}(\glL)}
\def\UM{\operatorname{U}(\glM)}
\def\UMN{\operatorname{U}(\glMN)}
\def\UN{\operatorname{U}(\glN)}
\def\Up{\Upsilon}
\def\Upt{\widetilde{\Upsilon}}

\def\Vla{V_{\la\ts\lat}}
\def\VLa{V_{\ns\La\Lat\ts}}
\def\Vlm{V_{\ns\la\lat}^{\,\mu\ts\mut}}
\def\VOm{V_{\Om\Omt}}

\def\Wkk{W_{k\ts\kt}}
\def\Wkko{W_{k\ts\kt}^{\hskip1pt\circ}}
\def\Wll{W_{l\ts\lt}}
\def\Wllo{W_{l\ts\lt}^{\hskip.5pt\circ}}
\def\Wmm{W_{m\mt}}
\def\Wmmo{W_{m\mt}^{\ts\circ}}
\def\Wnn{W_{n\ts\nt}}
\def\Wnno{W_{n\ts\nt}^{\,\circ}}
\def\Wp{W^{\hskip.5pt\prime}}

\def\YL{\operatorname{Y}(\glL)}
\def\YMN{\operatorname{Y}(\glMN)}
\def\YN{\operatorname{Y}(\glN)}

\def\Zmm{Z_{m\mt}}
\def\Zp{Z^{\hskip.5pt\prime}}

\def\ZZ{{\mathbb Z}} 

\begin{document}
\title{Rational representations of Yangians\\
associated with skew Young diagrams\\}
\author{Maxim Nazarov}
\institute{Department of Mathematics, University of York,
York YO10 5DD, England\\\email{mln1@york.ac.uk\\ }}
\titlerunning{Rational representations of Yangians}
\authorrunning{Maxim Nazarov}
\maketitle


\vskip-15pt
\begin{abstract}
{\bf\hskip-10pt.}\hfill  
Consider the general linear group $GL_M$ over the complex field. 
\newline
The irreducible rational representations of the group $GL_M$
can be labeled~by the pairs of partitions 
$\mu=(\mu_{\ts1},\mu_{\ts2}\ts,\,\ldots)$
and
$\mut=(\mut_{\ts1},\mut_{\ts2}\ts,\,\ldots)$
such that the~total number of non-zero parts of $\mu$ and $\mut$
does not exceed $M$. Let $V_{\mu\ts\mut}$ be the irreducible 
representation corresponding to such a pair.

Regard the direct product
$GL_N\times GL_M$ as a subgroup of $GL_{N+M}\ts$.
Take any irreducible rational representation
$V_{\la\lat}$ of $GL_{N+M}\,$. 
The vector space
$$
\Vlm={\rm Hom}_{\,GL_M}(\ts V_{\mu\ts\mut}\ts\com V_{\la\lat}\ts)
$$
comes with a natural action of $GL_N\ts$. Put 
$n=\la_1-\mu_1+\la_2-\mu_2+\ldots$
and $\nt=\lat_1-\mut_1+\lat_2-\mut_2+\ldots\ $.
For any pair of standard Young tableaux 
$\Om$ and $\Omt$
of skew shapes $\lm$ and $\lmt$ respectively,
we give a realization of $\Vlm$
as a subspace in the tensor product $\Wnn$ of $n$ copies of defining
representation $\CC^N$ of $GL_N\ts$, and of $\nt$ copies
of the contragredient representation $(\CC^N)^\ast$.
This subspace is determined as the image of a certain linear operator
$\FOm$ on $\Wnn\,$. We introduce this operator
by an explicit multiplicative formula.

When $M=0$ and $\Vlm=V_{\la\lat}$ is an irreducible representation of
$GL_N$, we recover the known realization of $V_{\la\lat}$
as a certain subspace in the space of all traceless tensors in $\Wnn\,$.
Then the operator $\FOm$ may be regarded as the rational analogue
of the classical Young symmetrizer, corresponding to the
standard tableau $\Om$ of shape $\la\ts$.
That symmetrizer is a certain linear operator on
the tensor product of $n$ copies of $\CC^N$,
with the image equivalent to the irreducible
polynomial representation of $GL_N$, 
corresponding to the partition $\la\,$. Even when
$M=0$, our formula for the operator $\FOm$ is new.

Our results are applications of the representation theory
of the Yangian $\YN$ of the Lie algebra $\glN\ts$.
In particular, $\FOm$ is an intertwining operator between
certain representations of the algebra $\YN$ on $\Wnn\,$.
We also introduce the notion of a rational representation
of the Yangian $\YN\ts$. As a representation of $\YN\ts$,
the image of $\FOm$ is rational and irreducible.
\end{abstract}


\section{\hskip-3.5pt. Main results}


\textbf{1.1.}
For any positive integer $N$, consider the general linear group $GL_N$
over the complex field $\CC\ts$. Consider also the corresponding
Lie algebra $\glN\ts$. Let us choose the Borel subalgebra in $\glN$ 
consisting of the upper triangular matrices, and fix the basis of the 
diagonal matrix units $E_{11}\lcd E_{NN}$ in the corresponding Cartan 
subalgebra of $\glN\ts$.
Now take any irreducible {\it rational} representation
of the group $GL_N\ts$. This representation is finite-dimensional,
and its matrix elements in any basis are rational functions of the
standard coordinates on $GL_N\ts$. Moreover, denominators of
these rational functions are powers of the determinant.
This representation can also be regarded
as a $\glN$-module.
The highest weight of this $\glN$-module is a non-increasing
sequence of $N$ integers. Let us write this sequence of integers as
$$
(\ts\nu_{\ts1}\lcd\nu_{p}
\com0\lcd0\com
-\nut_{q}\lcd-\nut_{\ts1}\ts)
$$
where $\nu_{\ts1}\ge\ldots\ge\nu_{p}>0$ and
$\nut_{\ts1}\ge\ldots\ge\nut_{q}>0\ts$.
It is convenient to put $\nu_{\ts i}=0$ for every index $i>p\ts$,
and $\nut_{\ts i}=0$ for every $i>q\ts$.
Thus our irreducible representation of $GL_N$ is
indexed by the pair of partitions
$$
\nu=(\ts\nu_{\ts1},\nu_{\ts2},\,\ldots\,)
\quad\text{and}\quad
\nut=(\ts\nut_{\ts1},\nut_{\ts2},\,\ldots\,)
$$
such that $\nup_{\ts1}+\widetilde\nu^{\,\prime}_{\ts1}\le N$.
As usual, $\nup=(\nup_{\ts1},\nup_{\ts2},\,\ldots\,)$
is the partition conjugate to $\nu\ts$.
In particular, here $\nup_{\ts1}=p$ is the number of non-zero parts of $\nu$.

Let us denote by $V_{\nu\ts\nut}$
irreducible representation of $GL_N$
corresponding to the partitions $\nu$ and $\nut$. Put 
$$
n=\nu_{\ts1}+\nu_{\ts2}+\ldots
\quad\text{and}\quad\,
\nt=\nut_{\ts1}+\nut_{\ts2}+\ldots\,
$$
so that $\nu$ and $\nut$ are partitions of the numbers $n$ and $\nt$
respectively.
The irreducible representation $V_{\nu\ts\nut}$ occurs in the tensor product
\begin{equation}\label{mix}
\Wnn=(\CC^N)^\ast\ot\ldots\ot(\CC^N)^\ast\ot\CC^N\ot\ldots\ot\CC^N
\end{equation}
of $n$ copies of the defining representation $\CC^N$ of the group $GL_N$,
and of $\nt$ copies of the contragredient representation $(\CC^N)^\ast\ts$.
The vector space at the right hand side of (\ref{mix}) is called a
space of {\it mixed tensors\/}. 
If $\nut=(0\com0\ts,\ts\ldots\ts)\ts$ and $\nt=0\ts$, then
$V_{\nu\ts\nut}\ts$ is a {\it polynomial\/} representation of $GL_N\ts$;
it is then denoted by $V_\nu\ts$. 
If $\nu=(0\com0\ts,\ts\ldots\ts)\ts$ and $n=0\ts$, then
the representation $V_{\nu\ts\nut}\ts$ of $GL_N$ is contragredient to
$V_{\ts\nut}\,$; it is then denoted by $V_{\nut}^{\ts\ast}\ts$ accordingly. 

Consider the natural pairing  
$\langle\ ,\,\rangle$ between the vector spaces
$(\CC^N)^\ast$ and $\CC^N$. Take any two indices
$k\in\{1\lcd n\}$ and $l\in\{1\lcd\nt\ts\}\ts$.
By applying the pairing $\langle\ ,\,\rangle$ to a tensor $w\in\Wnn$
in the $l\ts$th factor $(\CC^N)^\ast$ and the $k\ts$th factor $\CC^N$,
we obtain a certain tensor $\widehat{w}\in W_{n-1,\nt-1}\ts$.
The tensor $w$ is called \textit{traceless}, 
if $\widehat{w}=0$ for all possible $k$ and $l\ts$.
Denote by $\Wnno$ the subspace in $\Wnn$ consisting of all
traceless tensors, this subspace is $GL_N\ts$-invariant.
For the above $n$ and $\nt$, the irreducible components of $\Wnn$
equivalent to $V_{\nu\ts\nut}\ts$ occur only in the subspace $\Wnno\ts$.
Moreover, any such component of $\Wnn$ has the~form 
\vglue-16pt
$$
(\ts V_{\nut}^{\ts\ast}\otimes V_\nu\ts)\,\cap\,\Wnno
$$
for certain embeddings of the irreducible representations $V_\nu$ and 
$V_{\nut}^{\ts\ast}$ to the tensor products of $n$ copies 
of the $\CC^N$ and of $\nt$ copies of $(\CC^N)^\ast$ respectively.

All these facts about irreducible rational
representations of the group $GL_N$
are well known, see for instance \cite{K}.
In the present article, we consider a particular decomposition 
of the subspace of traceless tensors $\Wnno\subset\Wnn$ into
a direct sum of irreducible components.
In the particular case $\nt=0\ts$, this is the decomposition 
of the tensor product of $n$ copies of $\CC^N$ by means of
Young symmetrizers as in \cite[Section IV.4]{W}\ts;
see also \cite[Subsection 2.1]{N2}. For arbitrary $\nt\ts$, we
give an explicit multiplicative formula for the projector in $\Wnn$
onto each irreducible component of $\Wnno\ts$.
These projectors may be regarded as generalizations of the Young
symmetrizers from polynomial to the rational representations of $GL_N\ts$.
These generalizations seem to be new. 

Our multiplicative formulas for the generalized Young symmetrizers in
$\Wnn$ can be extended even further. Take any two pairs of partitions
$$
\la=(\ts\la_1,\la_2\ts,\,\ldots\,)
\quad\text{and}\quad
\lat=(\ts\lat_1,\lat_2\ts,\,\ldots\,)\,,
$$
$$
\mu=(\ts\mu_1,\mu_2\ts,\,\ldots\,)
\quad\text{and}\quad
\mut=(\ts\mut_1,\mut_2\ts,\,\ldots\,)\,.
$$
Then choose any integer $M$ such that 
\begin{equation}\label{ass}
\lap_{\ts1}+\widetilde\la^{\,\prime}_{\ts1}\le N+M
\quad\text{and}\quad
\mup_{\ts1}+\widetilde\mu^{\,\prime}_{\ts1}\le M\,.
\end{equation}
Consider the irreducible representations $V_{\la\lat}\ts$ 
and $V_{\mu\ts\mut}\ts$ of the groups $GL_{N+M}$ and $GL_M$
respectively. The decomposition $\ts\CC^{\ts N+M}=\CC^N\!\op\CC^M$ provides
an embedding of the direct product $GL_N\times GL_M$ to $GL_{N+M}$.
The vector space
\begin{equation}\label{1.0}
\Vlm\ts=\ts{\rm Hom}_{\,GL_M}(\,V_{\mu\ts\mut}\ts\com V_{\la\lat}\,)
\end{equation}
comes with a natural action of the group $GL_N$\ts.
This action of $GL_N$ may be reducible.
In the special case when $M=0$ and $\mu=\mut=(0\com0\ts,\ts\ldots\ts)\ts$,
we have $\Vlm=V_{\la\lat}\,$; we assume that $GL_M=\{1\}$ when $M=0$. 

The vector space $\Vlm$ is non-zero, if and only if
\begin{equation}\label{con}
\la_i\ge\mu_i\ ,\,
\lap_i-\mu^{\,\prime}_i\le N
\quad\text{and}\quad
\lat_i\ge\mut_i\ ,\,
\lat^{\,\prime}_i-\mut^{\ts\,\prime}_i\le N
\end{equation}
for every index $i\,$; see \cite[Section I.5]{M}\ts.
In this article, we give distinguished embeddings of $\Vlm$
into the space of mixed tensors $\Wnn$ where
$$
n=\la_1-\mu_1+\la_2-\mu_2+\ldots
\quad\text{and}\quad
\nt=\lat_1-\mut_1+\lat_2-\mut_2+\ldots\,.
$$
These embeddings will be compatible with the action of the group $GL_N\ts$.
Under any of these embeddings, the vector space $\Vlm$ will be realized
as the image of a certain operator on $\Wnn$ commuting with the action
of $GL_N\ts$. We will give
a multiplicative formula for each of these operators.
Setting $\mu=\mut=(0\com0\ts,\ts\ldots\ts)$ and $M=0$
in these multiplicative formulas, we will obtain projectors
onto $GL_N\ts$-irreducible components of $\Wnno\ts$.
However, for arbitrary $\mu$ and $\mut$ the images of
$\Vlm$ under our embeddings may be not contained in the subspace
of traceless tensors $\Wnno\subset\Wnn\ts$.

When referring to theorems, propositions, lemmas and corollaries,
we will indicate the subsections where they respectively appear. 
There will be no more than one of each of these in any subsection,
so our referring system should cause no confusion. 
For example, the proposition and the lemma 
in the next subsection will be referred to as Proposition 1.2 and Lemma 1.2.


\smallskip\medskip\noindent\textbf{1.2.}
First consider the case $\lat=\mut=(0\com0\ts,\ts\ldots\ts)$
and $\nt=0\ts$. In this case the vector spaces (\ref{mix}) and
(\ref{1.0}) will be denoted simply
by $W_n$ and $V_\la^{\ts\mu}$.
Suppose that $\la_i\ge\mu_i$ for every index $i\ts$.
Consider the {\it skew Young diagram}
$$
\lm=\{\,(i\com j)\in\ZZ^2\ |\ i\ge1,\ \la_i\ge j>\mu_i\,\}\,.
$$
If $\mu=(0\com0\ts,\ts\ldots\ts)\ts$, this is the conventional Young diagram
of the partition $\la\ts$.
We will employ the standard
graphical representation \cite[Section I.1]{M}
of Young diagrams on the plane $\RR^2$ with two matrix style coordinates.
Here the first coordinate increases from top to bottom, while the second
coordinate increases from left to right. The element $(i\com j)\in\lm$
is represented by the unit box with the bottom right corner
at the point $(i\com j)\in\RR^2$. 

The set $\lm$ consists of $n$ elements.
A \textit{standard tableau} of shape $\lm$ is a bijection
$\Om:\lm\to\{1\lcd n\}$ such that $\Om(i\com j)<\Om(i+1\com j)$ and
$\Om(i\com j)<\Om(i\com j+1)$ for all possible $i$ and $j\ts$.
For every $k=1\lcd n$ put $c_k(\Om)=j-i$ if $k=\Om(i\com j)\ts$.
The difference $c_k(\Om)$ is the \textit{content} of the box
occupied by the number $k$ in the tableau $\Om\ts$.
Graphically, the tableau $\Om$ is represented by placing the
numbers $\Om(i\com j)$ into the corresponding boxes of $\lm$
on the plane $\RR^2$. 
By filling the boxes with the numbers
$1\lcd n$ by columns from the left to the right,
downwards in every column, we get
the \textit{column tableau} $\Om$ of shape $\lm\ts$.
This tableau is standard. Below on the left
we represent the column tableau $\Om$ for
$\la=(5,\ns3,\ns3,\ns3,\ns3,\ns0,\ns0,\ts\ldots)$ and
$\mu=(3,\ns3,\ns2,\ns0,\ns0,\ts\ldots)\ts$,
on the right we show the contents of the boxes of $\lm\ts$:

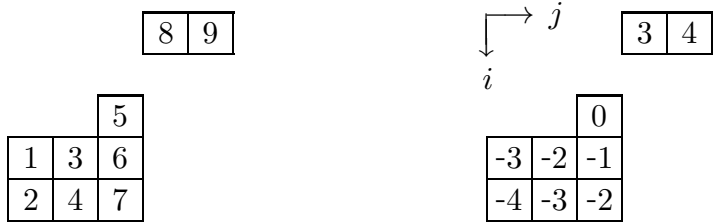
\bigskip
\vbox{
$$
\kern22.5pt\longrightarrow\,j\mw
$$
\vglue-25pt
$$
\kern80pt\vert
$$
\vglue-27pt
$$
\bigr\downarrow\mw\kern1pt
$$
\vglue-16pt
$$
i\mw
$$
\vglue-44pt
$$
\phantom{\bx}
\phantom{\bx}
\phantom{\bx}
{\bx}
{\bx}
\kern80pt
\phantom{\bx}
\phantom{\bx}
\phantom{\bx}
{\bx}
{\bx}
$$
\vglue-16.6pt
$$
\phantom{\bx}
\phantom{\bx}
\phantom{\bx}
\phantom{\bx}
\phantom{\bx}
\kern80pt
\phantom{\bx}
\phantom{\bx}
\phantom{\bx}
\phantom{\bx}
\phantom{\bx}
$$
\vglue-16.6pt
$$
\phantom{\bx}
\phantom{\bx}
{\bx}
\phantom{\bx}
\phantom{\bx}
\kern80pt
\phantom{\bx}
\phantom{\bx}
{\bx}
\phantom{\bx}
\phantom{\bx}
$$
\vglue-16.5pt
$$
{\bx}
{\bx}
{\bx}
\phantom{\bx}
\phantom{\bx}
\kern80pt
{\bx}
{\bx}
{\bx}
\phantom{\bx}
\phantom{\bx}
$$
\vglue-16.5pt
$$
{\bx}
{\bx}
{\bx}
\phantom{\bx}
\phantom{\bx}
\kern80pt
{\bx}
{\bx}
{\bx}
\phantom{\bx}
\phantom{\bx}
$$
\vglue-82.5pt
$$
\kern43pt8\kern9pt9\kern90pt\kern42pt3\kern9pt4
$$
\vglue-3.5pt
$$
5\kern146pt0
$$
\vglue-16.5pt
$$
1\kern9pt3\kern9pt6\kern116pt\mi3\kern6pt\mi2\kern6pt\mi1\kern26pt
$$
\vglue-16.5pt
$$
2\kern9pt4\kern9pt7\kern116pt\mi4\kern6pt\mi3\kern6pt\mi2\kern26pt
$$
}

\bigskip\medskip\noindent
In this example $n=9$, and the sequence of contents
$c_1(\Om)\lcd c_9(\Om)$ is
$$
(\ts\mi3\com\mi4\com\mi2\com\mi3\com0\com\mi1\com\mi2\com3\com4\ts)\ts.
$$

Now introduce complex
variables $t_1(\Om)\lcd t_n(\Om)$ with the constraints
\begin{equation}\label{1.0002}
t_k(\Om)=t_l(\Om)
\text{\ \ if $k$ and $l$ occur in the same column of $\Om$\ts.}
\end{equation}
So the number of independent variables among $t_1(\Om)\lcd t_n(\Om)$ 
equals the number of non-empty columns of the diagram $\lm\ts$.
Order lexicographically
the set of all pairs $(k\com l)$ with $1\le k<l\le n$.
Take the ordered product over this set,

\vskip-20pt
\begin{equation}\label{1.3}
\prod_{1\le k<l\le n}^{\longrightarrow}\ 
\left(1-\frac{P_{kl}}{\ts c_k(\Om)-c_{\ts l}(\Om)+t_k(\Om)-t_l(\Om)}\ts\right)
\end{equation}
\noindent
where $P_{kl}$ denotes the operator on the vector space $W_n$
exchanging the $k\ts$th and $l\ts$th tensor factors. Let us consider 
(\ref{1.3}) as a rational function of the constrained variables 
$t_1(\Om)\lcd t_n(\Om)$. The next result goes back to \cite{C1}.

\begin{Lemma}
The product {\rm(\ref{1.3})} is regular at\/
$t_1(\Om)=\ldots=t_n(\Om)\ts$. 
\end{Lemma}

\noindent
Note that the function (\ref{1.3}) depends only on the differences
$t_k(\Om)-t_l(\Om)\ts$. The value of {\rm(\ref{1.3})} at
$t_1(\Om)=\ldots=t_n(\Om)\ts$ will be denoted by $F_\Om\,$. 
For the proof of Lemma 1.2 see \cite[Section 2]{NT2}.
This proof provides an explicit multiplicative formula for
the operator $F_\Om$ on the vector space $W_n\ts$.

Let us denote by $V_\Om$ the image of the operator $F_\Om$ on $W_n\ts$.
The subspace $V_\Om\subset W_n$ is stable under the action of
$GL_N\ts$, see the last paragraph of this subsection. 
Proof of the following result was given in \cite[Subsection 4.6]{N2}.

\begin{Proposition}
The representations $V_\la^{\ts\mu}$ and $V_\Om$ of\/ $GL_N$ are equivalent. 
\end{Proposition}

\noindent
Note that the operator $F_\Om$ on the vector space $W_n$ does not depend
on $M$. It is well known that the dimension of the vector
space $V_\la^{\ts\mu}$ is the same for all integers $M$ such that
$\lap_1\le N+M$ and $\mup_1\le M$; see \cite[Section I.5]{M}.

The symmetric group $S_n$ acts on the tensor product $W_n$ of
$n$ copies of $\CC^N$ by permutations of the tensor factors. 
This action extends to the group ring $\CC S_n\ts$. 
The image of $\CC\ts S_n$ in the operator algebra 
$\End(W_n)$ coincides with
the commutant of the image of the group $GL_N\ts$;
see \cite[Section IV.4]{W}.
The operator $P_{kl}\in\End(W_n)$ in
(\ref{1.3}) corresponds to the transposition $(k\ts l)\in S_n\ts$.
The operator $F_\Om\in\End(W_n)$ corresponds to a certain element
$f_{\ts\Om}\in\CC S_n$ which does not depend on $N$. 
If $\mu=(0\com0\ts,\ts\ldots\ts)\ts$, then $f_{\ts\Om}$ is a diagonal
matrix element of the irreducible representation of the group $S_n$  
labeled by the partition $\la\,$; see Subsection 2.2 of the present article.
For an alternative definition of the element $f_{\ts\Om}\in\CC S_n$
in the case $\mu\neq(0\com0\ts,\ts\ldots\ts)\ts$, see Subsection 2.3;
in particular, see the equality (\ref{2.9}).


\smallskip\medskip\noindent\textbf{1.3.}
Let us now formulate our results for the general partitions
$\lat$ and $\mut\ts$.
We will assume that $\la_i\ge\mu_i$ and $\lat_i\ge\mut_i$
for every index $k\ts$, otherwise we would have $\Vlm=\{0\}$.
Our embedding of the vector space $\Vlm$ into the space of mixed tensors
(\ref{mix}) will depend on the choice of two standard tableaux $\Om$ and $\Omt$
of the skew shapes $\lm$ and $\lmt$ respectively.

Choose any basis $v_1\lcd v_N$ in the vector space $\CC^N$.
Let $v_1^\ast\lcd v_N^{\ts\ast}$ be the dual basis in $(\CC^N)^\ast$ so that
$\langle\ts v_a^\ast\com v_b\ts\rangle=\de_{ab}$ for $a\com b=1\lcd N$.
The vector

\vskip-20pt
\begin{equation}\label{1.44444444}
w_0\,=\,\sum_{a=1}^N\,v_a^\ast\ns\ot v_a\,\in\,(\CC^N)^\ast\ns\ot\ts\CC^N
\end{equation}

\newpage\noindent
does not depend on the choice of the basis $v_1\lcd v_N$ and is invariant
under the action of $GL_N$ on $(\CC^N)^\ast\ns\ot\,\CC^N$.
Then consider the linear operator
\begin{equation}\label{1.45}
u\ot v\,\mapsto\,\langle\ts u\com v\ts\rangle\cdot w_0
\end{equation}
in the vector space
$(\CC^N)^\ast\ns\ot\,\CC^N$, it commutes with the action of $GL_N\ts$.

The tableaux $\Om$ and $\Omt$ define the sequences of contents
$c_1(\Om)\lcd c_n(\Om)\ns$ and $c_1(\Omt)\lcd c_{\ts\nt}(\Omt)\ts$.
For any indices $k\in\{1\lcd n\}$ and $l\in\{1\lcd\nt\ts\}\ts$,
denote by $Q_{\ts l,\nt+k}$ the linear operator on $\Wnn\ts$,
acting as (\ref{1.45}) in the $l\ts$th and $(\nt+k)\ts$th
tensor factors of (\ref{mix}), and acting
as the identity in the remaining $n+\nt-2$ tensor factors.
Consider the ordered product
\begin{equation}\label{1.333}
\GOm\ =\ 
\prod_{1\le k\le n}^{\longrightarrow}\ \biggl(\ 
\prod_{1\le l\le\nt}^{\longrightarrow}\ 
\biggl(1-\frac{Q_{\ts\nt-l+1\ts,\ts\nt+k}}
{\ts c_k(\Om)+{c_{\ts l}(\Omt)}^{\phantom{\prime}}\ns\!+N+M}
\,\biggr)\!\biggr)\,.
\end{equation}
Note that
$$
c_k(\Om)+{c_{\ts l}(\Omt)}\ge2-\lap_1-\lat_1^{\ts\prime}\ge2-N-M\ts,
$$
so every denominator in (\ref{1.333}) is non-zero.
Also consider the ordered product
$$
\GOmp\ =\ 
\prod_{1\le k\le n}^{\longleftarrow}\ \biggl(\ 
\prod_{1\le l\le\nt}^{\longleftarrow}\ 
\biggl(1-\frac{Q_{\ts\nt-l+1\ts,\ts\nt+k}}
{\ts c_k(\Om)+{c_{\ts l}(\Omt)}^{\phantom{\prime}}\ns\!+N+M}
\,\biggr)\!\biggr)\,,
\hskip21pt
$$
where the factors corresponding to the indices $k=1\lcd n$ and $l=1\lcd\nt$
are arranged from right to left, as indicated by the reversed arrows.

The tableau $\Omt$ defines
the linear operator $F_\Omt$
in the tensor product $W_{\ts\nt}$ of $\nt$ copies of $\CC^N$,
see Subsection~1.2.
Identify the dual vector space $W_{\ts\nt}^{\ts\ast}$ with
the tensor product of $\nt$ copies of $(\CC^N)^\ast$, so that
\begin{equation}\label{den}
\Wnn=W_{\ts\nt}^{\ts\ast}\otimes W_n\,.
\end{equation}
Under the permutational action of the group $S_{\ts\nt}\ts$
in $W_{\ts\nt}\ts$, the operator $F_\Omt$ on $W_{\ts\nt}$ corresponds to 
a certain element $f_\Omt$ of the group ring $\CC S_{\ts\nt}\ts$. 
The group $S_{\ts\nt}$ also acts by permutations of the tensor factors
in $W_{\ts\nt}^{\ts\ast}\ts$. The operator on $W_{\ts\nt}^{\ts\ast}$
corresponding to the element $f_\Omt\in\CC S_{\ts\nt}\ts$ will also be
denoted by $F_\Omt\ts$. However, we will always distinguish the vector
spaces $W_{\ts\nt}$ and $W_{\ts\nt}^{\ts\ast}\ts$. The linear operators
on $W_{\ts\nt}$ and $W_{\ts\nt}^{\ts\ast}$ reversing the order of the
$\nt$ tensor factors $\CC^N$ and $(\CC^N)^\ast$ respectively, will be
denoted by the same symbol $P_{\ts\nt}\ts$.
 
The tableau $\Om$ defines the linear operator $F_\Om$
in the tensor product $W_n\ts$, as in Subsection 1.2.
Introduce the linear operator on $\Wnn$
\begin{equation}\label{maindef}
\FOm\ts=\,\GOm\ts\cdot\ts((\ts P_{\ts\nt}\ts F_\Omt\,P_{\ts\nt})\ot F_\Om)\,.
\end{equation}

\begin{Lemma}
We have the equality of operators on  the vector space $\Wnn$
$$
\FOm\ts=\ts((\ts P_{\ts\nt}\ts F_\Omt\,P_{\ts\nt})\ot F_\Om)\ts\cdot\ts\GOmp\,.
$$
\end{Lemma}

We will prove Lemma 1.3 in Subsection 4.2.
Denote by $\VOm$
the image of the operator $\FOm$ on $\Wnn\ts$.
Note that unlike the subspace $V_\Om\ts\subset W_n\ts$, 
the subspace $\VOm\subset\Wnn$ may depend on the choice of $M$.
The subspace $\VOm\subset\Wnn$ is stable under the action of 
$GL_N\ts$. Here is a generalization of Proposition 1.2 
from polynomial to the rational representations of $GL_N$.

\begin{Proposition}
The representations $\Vlm$ and $\VOm$ of\/ $GL_N$ are equivalent. 
\end{Proposition}

The proof of Proposition 1.3 is given in the end if Subsection 5.4. 
Denote by $V^{\ts\ast}_\Omt$ 
the image of the operator $F_\Omt$ in
$W_{\ts\nt}^{\ts\ast}\ts$.
Due to Lemma 1.3, the image $\VOm$ of the operator $\FOm$
is contained in the subspace
$$
(\ts P_{\ts\nt}\cdot V^{\ts\ast}_\Omt\ts)\ot V_\Om\,\subset\,
W_{\ts\nt}^{\ts\ast}\otimes W_n\,.
$$
If $M=0$, then $\mu=\mut=(0\com0\ts,\ts\ldots\ts)\ts$. Then
the vector spaces $V_\Om$ and $P_{\ts\nt}\cdot V^{\ts\ast}_\Omt$
are equivalent respectively to $V_\la$ and $V_{\ts\lat}^{\ts\ast}$ as
representations of the group $GL_N$. We have the equality
\begin{equation}\label{1.4444}
\VOm=
((\ts P_{\ts\nt}\cdot V^{\ts\ast}_\Omt)\ot V_\Om)\ts\cap\ts\Wnno
\quad\text{if}\quad
M=0\,;
\end{equation}
see Subsection 3.4. If $M\neq0$,
the image $\VOm$ of the operator $\FOm$ may be not contained in
the subspace $\Wnno\subset\Wnn$ of traceless tensors.
Yet our proof of Proposition 1.3 is based on the equality (\ref{1.4444}).

If $M\neq0$, the operator $\FOm$
can also be defined via the equality (\ref{3.4444}),
cf.\ the equality (\ref{2.9}).
Both definitions of the operator $\FOm$ are motivated by the
representation theory of Yangians~\cite{MNO}, see the next subsection.


\smallskip\medskip\noindent\textbf{1.4.}
The vector space $\Vlm$
is irreducible under the natural action of
the subalgebra of $GL_M$\ts-invariants in the enveloping algebra
$\UMN\ts$. Denote this subalgebra by $\AMN\ts$,
it coincides with the centralizer of the subalgebra $\UM\subset\UMN\ts$.
Our Theorem 1.5 below describes the action of the algebra $\AMN$ on $\Vlm$ 
explicitly, using the \textit{Yangian} $\YN$ of the Lie algebra $\glN\ts$.   
The Yangian $\YN$ is a deformation of the
enveloping algebra of the polynomial current Lie algebra $\glN[x]$
in the class of Hopf algebras. 

The unital associative algebra $\YN$ has a 
family of generators $T_{ab}^{(i)}$ where $i=1,2,\ts\ldots\ts$ and
$a\com b=1\lcd N$. The defining relations for these generators
can be written by using the formal power series
\begin{equation}\label{1.31}
T_{ab}(x)=
\de_{ab}\cdot1+T_{ab}^{(1)}x^{-\ns1}+T_{ab}^{(2)}x^{-\ns2}+\,\ldots
\,\in\,\YN\,[[x^{-1}]]\,.
\end{equation}
Here $x$ is the formal parameter. Let $y$ be another formal parameter,  
then the defining relations in the associative algebra $\YN$
can be written as
\begin{equation}\label{1.32}
(x-y)\cdot[\ts T_{ab}(x)\ts,T_{cd}(y)\ts]\ts=\;
T_{cb}(x)\ts T_{ad}(y)-T_{cb}(y)\ts T_{ad}(x)\,,
\end{equation}
where $a\com b\com c\com d=1\lcd N\ts$.
The square brackets in (\ref{1.32}) denote the usual commutator.
If $N=1$, the algebra $\YN$ is commutative.
Using the series (\ref{1.31}),
the coproduct $\De:\YN\to\YN\ot\YN$ is defined by
\begin{equation}\label{1.33}
\De\bigl(T_{ab}(x)\bigr)\,=\,\sum_{c=1}^N\ T_{ac}(x)\ot T_{cb}(x)\,;
\end{equation}
the tensor product on the right hand side of the equality (\ref{1.33})
is taken over the subalgebra $\CC[[x^{-1}]]\subset\YN\,[[x^{-1}]]\ts$.
The counit homomorphism $\varepsilon:\YN\to\CC$ is determined by
the assignment $\,\varepsilon:\,T_{ab}(x)\ts\mapsto\ts\de_{ab}\cdot1$.

The antipode $\operatorname{S}$ on $\YN$ can be defined by using the element
\begin{equation}\label{1.71}
T(x)\,\ts=\sum_{a,b=1}^N\, E_{ab}\ot
T_{ab}(x)\in\End(\CC^N)\ot\YN\,[[x^{-1}]]\,,
\end{equation}
where the matrix units $E_{ab}$ are regarded as basis elements of the
algebra $\End(\CC^N)\ts$. The formal power series (\ref{1.71}) in
$x^{-1}$ is invertible, because its leading term is $1\ot1$.
The anti-automorphism $\operatorname{S}$ is defined by the assignment
\begin{equation}\label{antip}
\id\ot\operatorname{S}\ts:\ts T(x)\mapsto T(x)^{-1}\,.
\end{equation}
We will also use an automorphism $\om_N$ of the algebra $\YN$ defined by
the assignment

\vskip-20pt
\begin{equation}\label{1.51}
\id\ot\om_N\ts:\ts T(x)\mapsto T(-x)^{-1}\,.
\end{equation}
The automorphism $\om_N$ is clearly involutive.
For references and more details on the definition of the Yangian $\YN$
see \cite[Section 1]{MNO}\ts.

Take any formal power series $g(x)\in\CC[[x^{-1}]]$ with the leading
term $1$. The assignment

\vskip-20pt
\begin{equation}\label{1.61}
T_{ab}(x)\mapsto\,g(x)\cdot T_{ab}(x)
\end{equation}
defines an automorphism of the algebra $\YN\ts$,
see (\ref{1.31}) and (\ref{1.32}).

The Yangian $\YN$ contains the enveloping
algebra $\UN$ as a Hopf subalgebra. Regard 
the matrix units $E_{ab}\in\glN$ as generators of
the algebra $\UN\ts$. Due to the defining relations (\ref{1.32}),
an embedding of associative algebras $\UN\to\YN$ 
can be determined by the assignment
\begin{equation}\label{4.4}
E_{ab}\mapsto-\ts T_{ba}^{(1)}.
\end{equation}
This particular embedding will be used in Section 4, see Proposition 4.1.
This is a Hopf algebra embedding, because by (\ref{1.33})
we have the equality
\begin{equation}\label{Hopfemb}
\De\ts\bigl(\ts T_{ab}^{(1)}\bigr)\ts=\,
T_{ab}^{(1)}\ns\ot\ts1+1\ot\ts T_{ab}^{(1)}\,.
\end{equation}
The relations (\ref{1.32}) also show that the assignment
\begin{equation}\label{1.52}
\pi_N:\ts T_{ab}(x)\,\mapsto\,\de_{ab}\cdot1-E_{ba}\,x^{-1}
\end{equation}
defines a homomorphism $\pi_N:\YN\to\UN$
of associative algebras. By definition,
the homomorphism $\pi_N$ acts on the subalgebra
$\UN\subset\YN$ as the identity. However,
$\pi_N$ is not a Hopf algebra homomorphism.

The defining relations (\ref{1.32}) show that for any $z\in\CC\,$,
the assignment
\begin{equation}\label{tau}
\tau_z:\,T_{ab}(x)\ts\mapsto\,T_{ab}(x-z)
\quad\textrm{for all}\quad
a\com b=1\lcd N
\end{equation}
defines an automorphism $\tau_z$ of the algebra $\YN\ts$. Here the formal
power series $T_{ab}(x-z)$ in $(x-z)^{-1}$ should be re-expanded in $x^{-1}$.
Note that the automorphism $\tau_z$ acts on the subalgebra
$\UN\subset\YN$ as the identity.

By pulling the action of 
the enveloping algebra $\UN$ in the defining $\glN$-module $\CC^N\ts$
back through the homomorphism
$$
\pi_N\circ\,\tau_z:\YN\to\UN\,,
$$
we obtain a module over the algebra $\YN\ts$, 
called an \textit{evaluation module\/}. 
To indicate the dependence on the
parameter $z\ts$, let us denote this $\YN\ts$-module by $V(z)\ts$.
As a $\glN$-module, this is still the defining module $\CC^N$;
here we use the embedding (\ref{4.4}). Finally, by
using the antipode $\operatorname{S}$ of $\YN\ts$, 
define the $\YN\ts$-module $V(z)^\ast$ dual to $V(z)\ts$. 
As a $\glN$-module, $V(z)^\ast$ is the contragredient module 
$(\CC^N)^\ast$. Here we again use the embedding (\ref{4.4}). 

The operator $\FOm$ on the space
$\Wnn$ admits the following interpretation in terms of the tensor products
of evaluation modules over the algebra $\YN\ts$. Here we will use
the comultiplication (\ref{1.33}) on $\YN\ts$. Let 
\begin{equation}\label{Pnn}
P_{\ts n+\nt}\ts:\ts
W_n\ot W_{\ts\nt}^{\ts\ast}\,\to\,
W_{\ts\nt}^{\ts\ast}\ot W_n
\end{equation}
be the linear
operator reversing the order of all the $\ts n+\nt\ts$ tensor factors. 

\begin{Proposition}
The operator\/ $\FOm\,P_{n+\nt}\ts$ is a\/ $\YN$-intertwiner
$$
V(c_n(\Om))\ot\ldots\ot V(c_1(\Om))\ot
V(-c_1(\Omt)\ns-\!M)^\ast\ot\ldots\ot V(-c_{\ts\nt}(\Omt)\ns-\!M)^\ast
$$
\vglue-10pt
$$
\Big\downarrow
$$
\vglue-14pt
$$
V(-c_{\ts\nt}(\Omt)\ns-\!M)^\ast\ot\ldots\ot V(-c_1(\Omt)\ns-\!M)^\ast\ot
V(c_1(\Om))\ot\ldots\ot V(c_n(\Om))\,.
$$
\end{Proposition}

\noindent
Thus the image $\VOm$ of the operator $\FOm$ on $\Wnn$
is a submodule of the target $\YN\ts$-module in Proposition 1.4.
If $\,\nt=0\ts$, this interpretation of $\FOm=F_\Om$ and $\VOm=V_\Om$
is due to Cherednik~\cite{C2}. In full generality,
the proof of Proposition 1.4 is given is Subsection 4.2
of the present article.

If $M=0$, then
$\mu=\mut=(0\com0\ts,\ts\ldots\ts)\ts$. 
The image $\VOm$ of the operator $\FOm$ is then equivalent to
$V_{\la\lat}$ as a representation of $GL_N$.
Proposition~1.4 then turns $V_{\la\lat}$ into $\YN\ts$-module. 
Note that the resulting $\YN\ts$-module can also be obtained from the
$\glN\ts$-module $V_{\la\lat}$ by pulling back through the
homomorphism $\pi_N\ts$, and then through a certain automorphism of $\YN$
of the form (\ref{1.61})\ts;
this is a special case of Theorem~1.5 below.


\smallskip\medskip\noindent\textbf{1.5.}
Olshanski \cite{O1} defined a homomorphism from the algebra $\YN$ 
to the subalgebra $\AMN$ of $GL_M$-invariants in $\UMN$,
for each non-negative integer $M$. Along with the centre
of the algebra $\UMN\ts$, the image of this homomorphism generates
the algebra $\AMN\ts$. We will use the following version of this homomorphism,
it will be denoted by $\pi_{NM}\ts$.

Let the indices $a\com b$ range over the set $\{1\lcd N+M\}\ts$.
Consider the basis of the matrix units $E_{ab}$
in the Lie algebra $\glMN\ts$. We assume that the subalgebras $\glN$ and
$\glM$ in $\glMN$ are spanned by elements $E_{ab}$ where
$$
1\le a\com b\le N
\quad\text{and}\quad
N+1\le a\com b\le N+M
$$
respectively. The subalgebra in the Yangian $\YMN$
generated by $T_{ab}^{(i)}$ where $1\le a\com b\le N\,$,
by definition coincides with the Yangian $\YN$.
Let us denote by $\ph_M$ this natural embedding $\YN\to\YMN$. 
Consider also the involutive automorphism $\om_{N+M}$
of the algebra $\YMN\ts$, see (\ref{1.51}).
The image of the homomorphism
$$
\pi_{N+M}\circ\ts\om_{N+M}\circ\ts\ph_M:\ts\YN\to\UMN
$$
belongs to the subalgebra $\AMN\subset\UMN\ts$.
Moreover, this image along with the centre
of the algebra $\UMN\ts$, generates the subalgebra $\AMN\ts$.
For the detailed proofs of these two assertions, see \cite[Section 2]{MO}.
In the present article, we use the homomorphism $\YN\to\UMN$
\begin{equation}\label{1.69}
\pi_{NM}=\,\pi_{N+M}\circ\ts\om_{N+M}\circ\ts\ph_M\circ\ts\om_N\ts.
\end{equation}
Note that when $M=0$, the homomorphism (\ref{1.69}) coincides with $\pi_N\ts$. 
Put
\begin{equation}\label{1.62}
g_\mu(x)\ =\ \prod_{i=1}^\infty\ 
\frac
{(\ts x-\mu_i+i\ts)(\ts x+i-1\ts)}
{(\ts x-\mu_i+i-1\ts)(\ts x+i\ts)}
\,\,.
\end{equation}
Here in the infinite product over $i\ts$, only finitely many
factors differ from $1$. Hence $g_\mu(x)$ is a rational function 
of $x$. We have $g_\mu(\infty)=1$, therefore $g_\mu(x)$
expands as a power series in $x^{-1}$ with leading term $1$. 

Let us keep to the assumptions \ns(\ref{ass}) on the integer $M$.
The~vector~space $\Vlm$ comes
with a natural action of the algebra $\AMN\ts$. Let us regard $\Vlm$ as a 
$\YN\ts$-module, by using the composition of the homomorphism
$\pi_{NM}:\YN\to\AMN$ with the automorphism of $\YN$ defined by
(\ref{1.61}) where 
$$
g(x)\,=\,g_\mu(x)\,\ts
g_{\ts\lat}(\ts-\ts x\ns-\!M)
\ts/\ts
g_{\ts\mut\ts}(\ts-\ts x\ns-\!M)\,.
$$
By Proposition~1.4, the image $\VOm$ of the operator
$\FOm$ can also be regarded as a $\YN\ts$-module.
Our main result is the following theorem.

\begin{Theorem}
The\/ $\YN$-modules\/ $\Vlm$ and\/ $\VOm$ are equivalent.
\end{Theorem}

\noindent
When $\nt=0\ts$, this theorem goes back to \cite[Theorem 2.6]{C2}.
However, most of the results in \cite{C2} have been given without proofs.
We prove Theorem~1.5 in Subsections 5.2 to 5.4 of the present article.
The algebra~$\AMN$ acts on the vector space $\Vlm$ irreducibly
\cite[Theorem 9.1.12]{D}.
The central elements of $\UMN$ act on $\Vlm$
as scalar operators. So Theorem 1.5 has a corollary.

\begin{Corollary}
The $\YN\ts$-module $\VOm$ is irreducible.
\end{Corollary}

The embedding $\UN\to\YN$ defined by (\ref{4.4}) 
provides an action of $\UN$ on the $\YN\ts$-module $\Vlm\ts$.
The vector space $\Vlm$ also comes with a natural action
of $\UN\subset\AMN\ts$. This natural action of $\UN$ in
$\Vlm$ coincides with its action as a subalgebra in $\YN\ts$,
see Subsection~4.3.

The subspace $\VOm\subset\Wnn$
is stable under the action of the group $GL_N\ts$.
Hence the enveloping algebra $\UN$ acts naturally on the vector space
$\VOm\ts$. This natural action of $\UN$ on $\VOm$ coincides with its action
as a subalgebra in $\YN\ts$, see again Subsection 4.3.

The natural action of $\UN$ on 
the vector space $\Vlm$ may be reducible.
Using Theorem 1.5 and its Corollary 1.5,
we can identify the vector space $\Vlm$
with the subspace $\VOm\subset\Wnn$ uniquely, up to multiplication
in $\VOm$ by a non-zero complex number.
Theorem 1.5 can be regarded as sharpening of Proposition 1.3.
Moreover, we will obtain Proposition 1.3 in the course of
the proof of Theorem 1.5. In our proof, we will use Proposition 3.6.


\smallskip\medskip\noindent\textbf{1.6.}
In priniciple, the analysis of any rational representation 
of the group $GL_N$ may be reduced to a polynomial representation,
cf.\ \cite[Section IV.5]{W}. In particular,
the irreducible rational
representation $V_{\nu\ts\nut}$ of $GL_N$ as defined in Subsection~1.1,
is equivalent to the tensor product of the one-dimensional
representation where

\vskip-16pt
$$
X\ts\mapsto\,(\ts\det X)^{-\nut_1}
$$
for any $\ts X\in GL_N\ts$, and of the polynomial representation $V_\zeta$ 
of $GL_N$ where
$$
\zeta_{\ts i}\ts=\,
\left\{
\begin{array}{ll}
\nu_{i}-\nut_{\ts N-i+1}+\nut_{\ts1}
&\ \ \textrm{if}\ \ \ i\le N,\\[2pt]
0
&\ \ \textrm{if}\ \ \ i>N.
\end{array}
\right.
$$
However, realization of the representation $V_{\nu\ts\nut}$
in the space of mixed tensors (\ref{mix}) is preferable,
at least from the combinatorial point of view \cite{KW,S}\ts,
and when $N\to\infty$ \cite{O2,VK}.
This makes our explicit formula
(\ref{maindef}) for the operator $\FOm$ on (\ref{mix})
useful, already 
when $M=0$ and $\mu=\mut=(0\com0\ts,\ts\ldots\ts)\ts$.

Borrowing terminology from the representation theory of the group
$GL_N$ let us call an irreducible
module over the Hopf algebra $\YN$ \textit{polynomial\/}, 
if it is equivalent to a submodule in the tensor product of
evaluation modules $V(z_1)\ot\ldots\ot V(z_n)$ for some integer $n\ge0$ and
some $z_1\lcd z_n\in\CC\ts$. 
Let us call an irreducible $\YN\ts$-module \textit{rational\/}, 
if it is equivalent to a submodule in a tensor product of 
any number of evaluation
and dual evaluation $\YN\ts$-modules. The tensor product here may be
taken in any suitable order, since the Hopf algebra $\YN$ is not
cocommutative for $N>1$. The irreducible $\YN\ts$-module $\VOm$ is then
rational by definition. The $\YN\ts$-module $\Vlm$ as figured in
Theorem 1.5 is then also rational due to Corollary 4.4.

Any irreducible rational
$\YN\ts$-module $V$ can be made polynomial by pulling $V$ back through
an automorphism of $\YN$ of the form (\ref{1.61}) for a suitable series 
$g(x)\in\CC[[x^{-1}]]$ with leading term $1$, cf. \cite[Theorem 2.16]{CP}.
Note that due to the definition (\ref{1.33}),
this pullback of $V$ coincides with the tensor product of $V$
with the one-dimensional $\YN\ts$-module such that
\begin{equation}\label{onedimrep}
T_{ab}(x)\,\mapsto\,\de_{ab}\cdot g(x)\,.
\end{equation}

Let us show how to make the $\YN\ts$-module $\VOm$ polynomial, 
according to the general assertion above.
The standard tableaux $\Om$ and $\Omt$ are of skew shapes
$\lm$ and $\lmt$ respectively. Suppose $\VOm\neq\{0\}$,
so that the conditions (\ref{con}) are satisfied.
Put $r=\lat_1\ts$. Define 
two partitions $\xi$ and $\eta$ by
$$
\xi_{\ts i}\ts=\,
\left\{
\begin{array}{ll}
\la_{i}-\lat_{\ts N+M-i+1}+r
&\ \ \textrm{if}\ \ \ i\le N+M,\\[2pt]
0
&\ \ \textrm{if}\ \ \ i>N+M
\end{array}
\right.
$$
and
$$
\eta_{\ts i}\ts=\,
\left\{
\begin{array}{ll}
\mu_{i}-\mut_{\ts M-i+1}+r
&\ \ \textrm{if}\ \ \ i\le M,\\[2pt]
0
&\ \ \textrm{if}\ \ \ i>M.
\end{array}
\right.
$$
Then we have
$$
\xi_{\ts i}\ge\eta_{\ts i}\ ,\,
\xi^{\,\prime}_{\ts i}-\eta^{\,\prime}_{\ts i}\le N
$$
for every index $i$ due to (\ref{con}).
In particular, the skew Young diagram $\xi/\eta$ is well defined.
Choose any standard tableau $\Gamma$ of shape $\xi/\eta\ts$.
Consider the corresponding $\YN\ts$-module $V_\Gamma\ts$,
it is polynomial by definition. Denote by $V_\Gamma(-r)$ the pullback
of $V_\Gamma$ relative to the automorphism $\tau_{-r}$ of $\YN\ts$.
The $\YN\ts$-module $V_\Gamma(-r)$ is also polynomial by definition.
The proof of the next proposition will be given in Subsection 4.5\ts;
cf.\ \cite[Section A]{KW}.

\begin{Proposition}
The\/ 
$\YN$-module\/ $\VOm$ is equivalent to the
pullback of the\/ $\YN\ts$-module\/
$V_\Gamma(-r)$ relative to the automorphism\/
{\rm(\ref{1.61})}, where
$$
g(x)\,=\,g_{\ts\lat}(\ts-\ts x\ns-\!M)
\,\cdot\,
\frac{\,x+M+r}{x+M}\ .
$$
\end{Proposition}

\noindent
Together with Theorem 1.5, this proposition provides another realization
of the rational $\YN\ts$-module $\Vlm$.
However, the realization of $\Vlm$
in the space (\ref{mix}) of mixed tensors provided directly by Theorem 1.5
is preferable, as it was for the rational representation
$V_{\nu\ts\nut}$ of the group $GL_N$. 

In the particular case when $\nt=0\ts$,
the $\YN\ts$-module $\VOm$ is already polynomial.
Proof of Theorem 1.5 in this case was given in \cite[\ns Section~4]{N2}. 
The present article is a sequel to \cite{N2}, but can be read independently.


\section{\hskip-3.5pt. Fusion procedure}


\textbf{2.1.}
Let us begin this section with recalling a few classical facts
about the irreducible representations of the symmetric group  $S_l$ over the
field $\CC$.  These representations are labeled by partitions $\la$ of $l\ts$.
We will identify partitions with their Young diagrams.
Denote by $U_\la$ the irreducible representation of $S_l$
corresponding to the partition $\la\ts$. We will also regard
representations of the group $S_l$ as modules over the group ring $\CC S_l$.
Fix the chain of subgroups $S_1\subset S_2\subset\ldots\subset S_l$
with the standard embeddings.

There is a decomposition of the vector space $U_\la$ into a direct sum of
one-dimensional subspaces, labeled by the {standard tableaux} of shape $\la$.
The one-dimensional subspace $U_\La\subset U_\la$ corresponding to
a standard tableau $\La$ is defined as follows. For any
$m\in\{1\lcd l-1\}$ take the tableau obtained from $\La$ by
removing the numbers $m+1\lcd l\ts$. Let the Young diagram $\mu$
be the shape of the resulting tableau.
Then the subspace $U_\La$ is contained in an irreducible
$\CC S_m$-submodule of $U_\la$ corresponding to $\mu\ts$. Any basis of
$U_\la$ formed by vectors $u_\La\in U_\La$ is called a 
\textit{Young basis}. Fix an $S_l$-invariant inner product
$(\,\,,\,)$ on $U_\la\ts$. All the subspaces
$U_\La\subset U_\la$ are then pairwise orthogonal. We choose the vectors
$u_\La\in U_\La$ so that $(\ts u_\La,u_\La\ts)=1$.

For any standard tableau $\La$ of shape $\la\ts$,
we will work with the diagonal matrix element of the
representation $U_\la$ corresponding to the vector $u_\La$\,,
\begin{equation}\label{2.0}
f_\La=
\sum_{s\in S_l}\,
(\,u_\La\com\ts s\ns\cdot\ns u_\La\,)
\,s\,\in\,\CC S_l\,.
\end{equation}
Note the equality
\begin{equation}\label{2.01}
f_\La^{\ts2}\,=\ts f_\La\cdot l\ts!\ts/\dim U_\la\ts.
\end{equation}

For any $k=1\lcd l-1$ let $s_k\in S_l$ be the transposition of $k$ and $k+1$.
We will use the description \cite{Y2} of the action of the generators 
$s_1\lcd s_{l-1}$ of the group $S_l$ on
the vectors of the Young basis.
Fix any standard tableau $\La\ts$. For every
$k=1\lcd l$ let $c_k(\La)$ be the {content}
of the box occupied by $k$ in $\La\ts$. That is, $c_k(\La)=j-i$ if
$k=\La(i\com j)\ts$.
Consider the tableau $s_k\La$ obtained from $\La$ by exchanging the numbers
$k$ and $k+1$. The tableau $s_k\La$ may be non-standard, this
happens exactly when $k$ and $k+1$ stand (next to each other) in the same 
row or column of $\La\ts$. But $c_k(\La)\neq c_{k+1}(\La)$ always, put
$h=(\ts c_{k+1}(\La)-c_k(\La))^{\ts-1}$. If $s_k\La$ is non-standard,
then $h=1$ or $h=-1$.

So far the vector $u_\La$ has been determined
up to a multiplier $z\in\CC$ with $|z|=1$. According to \cite{Y2}, all the 
vectors of the Young basis can be further normalized so that
for any standard tableau $\La$ and $k=1\lcd l-1$
\begin{equation}\label{2.3}
s_k\cdot u_\La=
\left\{
\begin{array}{ll}
h\ts u_\La+\sqrt{1-h^2}\, u_{s_k\La}
&\ \ \textrm{if}\ s_k\La\ \textrm{is standard,}\\[2pt]
h\ts u_\La
&\ \ \textrm{otherwise.}
\end{array}
\right.
\end{equation} 
This normalization fixes all the vectors of the Young basis
up to a common multiplier $z\in\CC$ with $|z|=1$. If the tableau $s_k\La$
is standard, then by (\ref{2.3}) 
\begin{equation}\label{2.55555}
\sqrt{1-h^2}\ts\,u_{s_k\La}=(s_k-h)\,u_\La\,.
\end{equation}
Then by the definition (\ref{2.0}) we have the identity
\begin{equation}\label{2.4}
(s_k-h)\,f_\La=f_{s_k\La}\,(s_k+h)\ts.
\end{equation} 
Note that for any two standard tableaux $\La$ and $\Lap$ there exists a
sequence of transpositions $s_{k_1}\lcd s_{k_n}$ such that
$\Lap=s_{k_n}\ldots s_{k_1}\La$ and the tableau $s_{k_m}\ldots s_{k_1}\La$
is standard for each $m=1\lcd n-1$. Then the identities (\ref{2.4})
provide a relation between the matrix elements $f_\La$ and $f_\Lap$
for any $\La\neq\Lap$.

We end the present subsection with the following useful identity.
Define the function $g_\la(x)\in\CC(x)$ via (\ref{1.61}).
Then by \cite[Lemma 4.4]{N2} we have
\begin{equation}\label{fact}
g_\la(x)\ =\ 
\prod_{k=1}^l\ 
\biggl(
1-\frac1{(\ts x-c_k(\La))^2}
\biggr)^{\!-1}\ts.
\end{equation}


\medskip\noindent\textbf{2.2.}
The matrix element $f_\La\in\CC S_l$
can also be obtained by the so called \textit{fusion procedure}.
For every two distinct indices
$i\com j\in\{1\lcd l\ts\}$ introduce the rational
function of $x\com y\in\CC$
\begin{equation}\label{2.45}
f_{ij}(x\com y)=1-\frac{(\ts i\ts j\ts)}{x-y}\ ,
\end{equation}
valued in $\CC S_l$\,; 
here $(\ts i\ts j\ts)\in S_l$ is the transposition of $i$ and $j\ts$.
As direct calculation shows, these rational functions satisfy the relation
\begin{equation}\label{2.5}
f_{ij}(x\com y)\,f_{ik}(x\com z)\,f_{jk}(y\com z)=
f_{jk}(y\com z)\,f_{ik}(x\com z)\,f_{ij}(x\com y)
\end{equation}
for any three pairwise distinct indices $i\com j\com k\in\{1\lcd l\ts\}\ts$. 
Note that
\begin{equation}\label{2.55}
f_{\ts ij}(x\com y)\,f_{\ts ji}(y\com x)=1-\frac1{(x-y)^{\ts2}}\,.
\end{equation}

Now take $l$ complex variables $x_1\lcd x_l$. Order lexicographically
the set of all pairs $(i\com j)$ with $1\le i<j\le l$. The
ordered product over this set,
\begin{equation}\label{2.6}
\prod_{1\le i<j\le l}^{\longrightarrow}\ 
f_{ij}(\ts x_i\com x_j\ts)
\end{equation}
is a rational function of $x_1\lcd x_l$ with values in $\CC S_l\ts$.
This rational function depends only on the differences $x_i-x_j\ts$.
Further, take the complex variables $t_1(\La)\lcd t_l(\La)$ 
constrained like in Subsection 1.2. That is,
\begin{equation}\label{1.22}
t_i(\La)=t_j(\La)
\text{\ \ if $i$ and $j$ occur in the same column of $\La$\ts.}
\end{equation}
So the number of independent variables amongst $t_1(\La)\lcd t_l(\La)$ 
is $\la_1\ts$. Set
\begin{equation}\label{1.2}
x_k=c_k(\La)+t_k(\La)
\text{\ \ for each\ \ $k=1\lcd l$\ts.}
\end{equation}

\begin{Proposition}
Restriction of the rational function {\rm(\ref{2.6})} to {\rm(\ref{1.2})}
is regular at $t_1(\La)=\ldots=t_l(\La)$. The value of this restriction at
$t_1(\La)=\ldots=t_l(\La)$ coincides with the element $f_\La\in\CC S_l$.
\end{Proposition}

\noindent
This proposition goes back to \cite{C2}.
In its present form, it has been proved in \cite[Section 2.2]{N1}.
The proof actually provides an explicit multiplicative formula for the
element $f_\La\in\CC S_l\ts$, different from the classical 
formula \cite{Y1}.


\medskip\smallskip\noindent\textbf{2.3.}
We need a generalization of Proposition 2.2
to standard tableaux of skew shapes \cite{C2}.
Take any $m\in\{0\lcd l-1\}$\ts. 
Let $\Up$ be standard tableau obtained from $\La$ by removing
the boxes with numbers $m+1\lcd l$. Let $\mu$ be the shape of
$\Up$. Define a standard tableau $\Om$ of skew shape
$\lm$ by setting 
\begin{equation}\label{Om}
\Om(i\com j)=\La(i\com\ns j)-m
\ \ \text{for each}\ \ 
(i\com j)\in\lm\,.
\end{equation}
Every standard tableau $\Om$ of shape $\lm$
can be obtained from a suitable $\La$ in this way.

Put $n=l-m$. Denote by $\io_m$ the embedding of the symmetric group $S_n$
into $S_l$ as a subgroup preserving the subset $\{m+1\lcd l\}$\ts;
we extend the mapping $\io_m$ to $\CC S_n$ by linearity.
Denote by $S_{mn}$ the subgroup $S_m\times\io_m(S_n)$ in $S_l\ts$.
Introduce the linear mapping
\begin{equation}\label{2.85}
\th_m:\CC S_l\to\CC S_{mn}:s\mapsto
\left\{
\begin{array}{ll}
s&\ \textrm{if}\ s\in S_{mn}\ts;\\[2pt]
0&\ \textrm{otherwise.}
\end{array}
\right.
\end{equation}
By definition, the element $f_\La\in\CC S_l$ is divisible on the
left and on the right by $f_{\ts\Up}\in\CC S_m$. Hence there exists an
element $f_\Om\in\CC S_n$ such that
\begin{equation}\label{2.9}
\th_m(f_\La)=f_{\ts\Up}\cdot\io_m(f_\Om).
\end{equation}
The element $f_\Om\in\CC S_n$ does not depend on the choice of
standard tableau $\Up$ of the shape $\mu$\ts, because the boxes with
the numbers $1\lcd m$ have in $\La$ and $\Up$ the same contents\ts;
see (\ref{2.4}). The generalization of Proposition 2.2 from $\La$
to $\Om$ is based on the following simple observation.

\begin{Proposition}
The image under the map\/ $\th_m$ of the product\/
 {\rm(\ref{2.6})} equals
\begin{equation}\label{2.7}
\prod_{1\le i<j\le m}^{\longrightarrow}
f_{ij}(\ts x_i\com x_j\ts)
\hskip5pt\cdot\hskip-6pt
\prod_{m<i<j\le l}^{\longrightarrow} 
f_{ij}(\ts x_i\com x_j\ts)\,.
\end{equation}
\end{Proposition}

\noindent
Now consider the ordered product on the right hand side~of~(\ref{2.7})\ts,
\begin{equation}\label{2.10}
\prod_{m<i<j\le l}^{\longrightarrow} 
f_{ij}(\ts x_i\com x_j\ts)\,.
\end{equation}

\begin{Corollary}
Restriction of the rational function {\rm(\ref{2.10})} to {\rm(\ref{1.2})} 
is regular at $t_{m+1}(\La)=\ldots=t_l(\La)$. 
The value at $t_{m+1}(\La)=\ldots=t_l(\La)$ of this restriction
coincides with $\io_m(f_\Om)\in\CC S_l$.
\end{Corollary}

\noindent
For the proof of Proposition 2.3 and Corollary 2.3,
see \cite[Subsection 2.3]{N2}. 

Let us complete the current subsection with another simple observation.
For any two standard tableaux $\La$ and $\Lap$ of the same shape $\la\ts$,
consider the matrix element of the representation $U_\la$ corresponding
to the pair of vectors $u_\La$ and $u_{\La'}$ of the Young basis,
\begin{equation}\label{fll}
f_{\ts\La}^{\ts\La'}=
\sum_{s\in S_l}\,
(\,u_{\La'}\com\ts s\cdot\ns u_\La\,)
\,s\,\in\,\CC S_l\,.
\end{equation}
If $\La=\Lap\ts$, then we have the equality $f_{\ts\La}^{\ts\La'}=f_\La$
by the definition (\ref{2.0}).
Now take any $m\in\{0\lcd l-1\}$. Let $n=l-m$ as before.
Consider the element
$$
\th_m(\ts f_{\ts\La}^{\ts\La'})\in S_{mn}\,,
$$
see (\ref{2.85}).
Let $\Up$ and $\Up^{\ts\prime}$ be the standard tableaux obtained by removing
the boxes with numbers $m+1\lcd l$ from the tableaux $\La$ and
$\Lap$ respectively. 
 
\begin{Lemma}
We have\/ $\th_m(\ts f_{\ts\La}^{\ts\La'})=0$ unless\/ 
$\Up$ and\/ $\Up^{\ts\prime}$ are of the same shape.
\end{Lemma}

\begin{proof}
By its definition, the element $ f_{\ts\La}^{\ts\La'}\in\CC S_l$ is divisible
by $f_{\ts\Up}\in\CC S_m$ on the right, and by $f_{\ts\Up'}\in\CC S_m$ on
the left. The image $\th_m(\ts  f_{\ts\La}^{\ts\La'})\in\CC S_{mn}$ 
inherits these two divisibility properties.
If the tableaux $\Up$ and\/ $\Up^{\ts\prime}$ are not of the same shape,
the irreducible representations $U_\Up$ and $U_{\Up'}$ of $S_m$ are not 
equivalent, and $f_{\ts\Up'}\ts s\,f_{\ts\Up}=0$ for any $s\in S_m$\qed
\end{proof}


\smallskip\noindent\textbf{2.4.}
Here we collect a few results that
we need for our proof of Theorem~1.5.
Let the symmetric group $S_{l+1}$ act by permutations on $0\com1\lcd l\ts$.
Then the subgroup $\io_1(S_l)\subset S_{l+1}$ still acts by permutations
on $1\lcd l\ts$. 
Consider the rational functions (\ref{2.45}) with
pairwise distinct indices $i\com j\in\{0\com1\lcd l\ts\}$;
these functions take values in the group ring $\CC S_{\ts l+1}\ts$.
Take the matrix element $f_\La\in\CC S_l$ defined by (\ref{2.0}).
Consider the image of $f_\La$ under the
embedding $\io_{\ts1}:\ts\CC S_l\to\CC S_{\ts l+1}\ts$.
For the proof of the next result see \text{\cite[Section 2]{N1}.}

\begin{Proposition}
We have equality of rational functions in\/ $x\ts$,
valued in\/ $\CC S_{\ts l+1}$
$$
\ f_{01}(x\com c_1(\La))\ts\ldots f_{0l}(x\com c_{\ts l}(\La))
\cdot\io_{\ts1}(f_\La)
\ts=\ts
\biggl(1\ts-\sum_{k=1}^l\,\frac{(\ts0\ts k)}x\ts\biggr)
\cdot\io_{\ts1}(f_\La)\,.
$$
\end{Proposition}

Take any $m\in\{0\lcd l-1\}$ and put $n=l-m\ts$, as in Subsection 2.3.
Also fix any non-negative integer $\lt$. Then take 
any $\mt\in\{0\lcd\lt-1\}$ and put $\nt=\lt-\mt\ts$.
Let the symmetric group $S_{1\ns+\lt+l}$ act by permutations on
\begin{equation}\label{inda}
0\com1\lcd\lt\com\lt+\ns1\lcd\lt+l\,.
\end{equation}
Define a linear map
$\ga:\CC S_{l+1}\to\CC S_{l+1}$ as follows.
Let the indices $i$ and $j$ range over the sequence (\ref{inda}).
Then for $s\in S_{l+1}$ we set
$\ga(s)=s\,$ if
$$
\lt<i\le\lt+m
\ \Rightarrow\ 
\lt<s(i)
\ \quad\text{and}\ \quad
\nt<j\le\lt
\ \Rightarrow\ 
0<s^{-1}(j)\le\lt\,;
$$
otherwise $\ga(s)=0$.
Furher, let $h_1\lcd h_{\lt+l}$ be arbitrary complex numbers.
Consider the ordered product in the algebra $\CC S_{1+\lt+l}\,$,
\begin{equation}\label{X}
\prod_{1\le i\le\lt+l}^{\longrightarrow}\!
((\ts0\ts i\ts)+h_i)\,.
\end{equation}

\begin{Lemma}
The image of the product\/ {\rm(\ref{X})}
under the linear map\/ $\ga$ equals
$$
\prod_{1\le i\le\nt}^{\longrightarrow}\,
((\ts0\ts i\ts)+h_i)
\hskip10pt\cdot\hskip-6pt
\prod_{\lt+m<i\le\lt+l}^{\longrightarrow}\hskip-7pt
((\ts0\ts i\ts)+h_i)\,.
$$
\end{Lemma}

\begin{proof}
Let us expand the ordered product (\ref{X}) as a sum of the products
of transpositions $s=(\ts0\ts i_1)\ldots(\ts0\ts i_a)$ with the coefficients
from $\CC\ts$; the sum is taken over all subsequences
$i_1\lcd i_a$ in the sequence $1\lcd\lt+l\ts$. Take any subequence such that
$\lt<i_b\le\lt+m$ for some $b\le a$. Assume that the index $b$ is minimal
with this property. Then $s(i_b)\le\lt$ so that $\ga(s)=0$.

Further, take any subsequence such that
$\nt<i_c\le\lt$ for some $c\le a$. Assume that the index $c$ is maximal
with this property. Then $s^{-1}(i_c)=0$ or 
$s^{-1}(i_c)>l\ts$, so that $\ga(s)=0$
\qed
\end{proof}

\noindent
We will also need a reformulation of this lemma.
Let us define a linear map
$\gap\ns:\ts\CC S_{l+1}\to\CC S_{l+1}$ by setting
$\gap(s)=s\,$ if
$$
\nt<i\le\lt
\ \Rightarrow\ 
0<s(i)\le\lt
\ \quad\text{and}\ \quad
\lt<j\le\lt+m
\ \Rightarrow\ 
s^{-1}(j)>\lt\,;
$$
otherwise $\gap(s)=0$.
Consider the ordered product in the algebra $\CC S_{1+\lt+l}\,$,
\begin{equation}\label{Y}
\prod_{1\le i\le\lt+l}^{\longleftarrow}\!
((\ts0\ts i\ts)+h_i)
\end{equation}
where the factors corresponding to the indices $i=1\lcd\lt+l$
are arranged from right to left, as indicated by the reversed arrow.

\begin{Corollary}
The image of the product\/ {\rm(\ref{Y})}
under the linear map\/ $\gap$ equals
$$
\prod_{\lt+m<i\le\lt+l}^{\longleftarrow}\hskip-7pt
((\ts0\ts i\ts)+h_i)
\hskip8pt\cdot\,\ts
\prod_{1\le i\le\nt}^{\longleftarrow}\,
((\ts0\ts i\ts)+h_i)\,.
$$
\end{Corollary}

\noindent
This corollary is derived from Lemma 2.4 by using the
anti-automorphism of the group ring $\CC S_{1+\lt+l}\,$, 
such that $s\mapsto s^{-1}$ for every group element $s\ts$.


\section{\hskip-3.5pt. Traceless tensors}


\textbf{3.1.}
For any positive integer $L$ take the vector spaces $\CC^L$ and
$(\CC^L)^\ast$. Then for any two non-negative integers $l$ and $\lt$
consider the tensor product
\begin{equation}\label{3.00}
\Wll\ts=\ts(\CC^L)^\ast\ot\ldots\ot(\CC^L)^\ast\ot\CC^L\ot\ldots\ot\CC^L
\end{equation}
of $l$ copies of the defining representation $\CC^L$ of the group $GL_L\ts$,
and of $\lt$ copies of the contragredient representation $(\CC^L)^\ast\ts$.
Note that now our basic vector space is $\CC^L$ while
in Section 1 it was $\CC^N\ns$, cf.\ (\ref{mix}).
The space $\CC^N$ will reappear in Subsection 3.5.
Denote by $\Bll$ the commutant of the image of $GL_L$
in the algebra $\End(\Wll)\ts$.
The operator algebra $\Bll$ may be regarded as a subalgebra of the 
\textit{Brauer centralizer algebra\/} \cite{B}, 
see for instance \cite{K}.

In this section, let us label the $\ts \lt+l\ts$ tensor factors 
in (\ref{3.00}) by the indices
\begin{equation}\label{ind}
1\lcd\lt\com\lt+\ns1\lcd\lt+l
\end{equation}
respectively from the left to the right.
Let $Q$ be the linear operator on $(\CC^L)^\ast\ot\CC^L$
defined similarly to the operator (\ref{1.45}) on $(\CC^N)^\ast\ot\CC^N$.
For any indices $i\com j\in\{1\lcd\lt+l\ts\}$ such that 
$i\le\lt<j\ts$, denote by $Q_{ij}$ the operator on $\Wll$ acting as $Q$
on the $i\ts$th and $j\ts$th tensor factors,
and acting as the identity on the remaining 
$\lt\ns+l\ns-\!2\ts$ tensor factors.
Each of the operators $Q_{ij}$ belongs to the commutant $\Bll\,$.
Note the equality $Q_{ij}^{\ts2}=L\cdot Q_{ij}\,$.

The subgroup $\Sll=S_{\ts\lt}\times\io_{\ts\lt\ts}(S_l)$ of the
symmetric group $S_{\ts\lt+l}$ acts on $\Wll$ by
permutations of $\lt+l$ tensor factors. 
The image of $\Sll$
in the algebra $\End(\Wll)$ belongs to $\Bll\ts$.
For two distinct indices $i\com j\in\{1\lcd\lt+l\ts\}$ the 
transposition $(i\ts j)\in S_{\ts \lt+l}$ belongs to the subgroup $\Sll$
if $\ts i\com j\le\lt\ts$ or if $\ts i\com j>\lt\ts$.
Then denote by $P_{ij}$ the corresponding permutational operator on 
$\Wll\,$. The centralizer algebra $\Bll$ is generated by operators $P_{ij}$
and $Q_{ij}$ with all possible indices $i\com j\,$;
see \cite[Lemma 1.2]{K}.   

The irreducible modules over the associative algebra $\Bll\ts$ 
or equivalently, the irreducible representations of the general linear 
group $GL_L$ that occur in the tensor product $\Wll\,$, can be indexed
by the pairs of partitions $\la$~and~$\lat$ of respectively 
$l-k$ and $\lt-k\ts$, where $k=0\lcd\min(\ts l\com\lt\ts)\ts$ and 
$\lap_{\ts1}+\lat^{\ts\prime}_{\ts1}\le L\ts$;
see for instance \cite[Corollary 4.7]{S}.
In the present article, we consider only the irreducible
$\Bll\,$-modules such that $k=0$.
The corresponding irreducible representations of the group $GL_L$
occur in the subspace of traceless tensors
$
\Wllo\,\subset\,\Wll\ts\,.
$
The images of the groups $GL_L$ and $\Sll$
in the algebra $\End(\Wllo\ts)$
span the commutants of each other
\cite[Theorem 1.1]{K}. In other words,
when acting on the vector space $\Wllo\,\ts$, 
the groups $GL_L$ and $\Sll$ form a \textit{dual pair}.

By definition, all
the operators $Q_{ij}$ on $\Wll$ vanish on the subspace $\Wllo\,$.
Denote by $\Ill$ the two-sided ideal in $\Bll$ generated by all the operators
$Q_{ij}\ts$. The quotient algebra $\Bll\ts/\,\Ill$ can be identified 
with the image of the group ring
$\CC\Sll$ in the algebra $\End(\Wllo\ts)\ts$.
The image of the element $P_{ij}\in\Bll$ in the quotient algebra 
is then identified with the operator on the subspace 
$\Wllo\,\subset\,\Wll\ts$
corresponding to $(i\ts j)\in\Sll\ts$. 

\vskip-2pt
Take any pair of partitions $\la$ of $l$ and $\lat$ of $\lt\ts$,
such that $\lap_{\ts1}+\lat^{\ts\prime}_{\ts1}\le L\ts$.
Denote by $\Vla$ the corresponding irreducible rational representation 
of the group $GL_L\ts$, this representation occurs in the subspace 
$\Wllo\,\ts$.
Further, take the irreducible representations $U_\la$ of $S_l$ and $U_{\lat}$ 
of $S_{\ts\lt}\ts$. Consider the tensor product 
$U_{\lat}\ot U_\la$ as a representation of the group $\Sll\ts$.
This is the irreducible representation, corresponding to $\Vla$
under the duality between $GL_L\ts$ and $\Sll\ts$ on $\Wllo\,\ts$.
Regard $U_{\lat}\ot U_\la$ as a $\Bll\,$-module
by using the homomorphism
\begin{equation}\label{3.22}
\Bll\,\longrightarrow\,\Bll\ts/\,\Ill\,.
\end{equation}
Then $U_{\lat}\ot U_\la$ is the irreducible $\Bll\,$-module
indexed by the pair $\la\ts\com\lat\ts$.


\smallskip\medskip\noindent\textbf{3.2.}
Let us keep labeling the tensor factors 
in (\ref{3.00}) by indices (\ref{ind}).
For any two distinct indices $i\com j\in\{1\lcd\lt+l\}$ such that 
$\ts i\com j\le\lt\ts$ or $\ts i\com j>\lt\ts$ put
\begin{equation}\label{3.45}
R_{\ts ij}(x\com y)=1-\frac{P_{ij}}{x-y}\,.
\end{equation}
As in Subsection 2.2, here $x$ and $y$ are complex variables. 
The function (\ref{3.45}) of $x$ and $y$ 
is called the \textit{Yang R-matrix\/}.
It corresponds to the function (\ref{2.45}) which now
takes values in $\CC S_{\ts\lt+l}\,$.
By (\ref{2.5}),
for any pairwise distinct indices $i\com j\com k\in\{1\lcd\lt+l\}$
such that $i\com j\com k\le\lt$ or $i\com j\com k>\lt$ we have the relation
\begin{equation}\label{3.5}
R_{\ts ij}(x\com y)\,R_{\ts ik}(x\com z)\,R_{\ts jk}(y\com z)=
R_{\ts jk}(y\com z)\,R_{\ts ik}(x\com z)\,R_{\ts ij}(x\com y)\,,
\end{equation}
it is called the \textit{Yang-Baxter relation\/}. Note that by (\ref{2.55}),
\begin{equation}\label{3.55}
R_{\ts ij}(x\com y)\,R_{\ts ji}(y\com x)=1-\frac1{(x-y)^{\ts2}}\,.
\end{equation}

For any two indices $i\com j\in\{1\lcd\lt+l\ts\}$ such that 
$i\le\lt<j\ts$, put
\begin{equation}\label{3.555}
\hskip40pt
\Rt_{ij}(x\com y)=1-\frac{Q_{ij}}{x-y}
\textrm{\hskip10pt and\hskip10pt}
\Rb_{ij}(x\com y)=1-\frac{Q_{ij}}{x+y+L}\,\ts;
\end{equation}
then we have

\vskip-20pt
\begin{equation}\label{3.6}
\Rt_{ij}(-x\com y)\ \Rb_{ij}(x\com y)=1\,.
\end{equation}

The linear operator $Q$ on $(\CC^L)^\ast\ot\CC^L$
is obtained from the permutation operator on $\CC^L\ot\CC^L$
by conjugation in the first tensor factor. Using this observation
along with the relation (\ref{2.5}), we get the relation
\begin{equation}\label{3.7}
\ \qquad
\Rt_{\ts ik}(x\com z)\,\Rt_{\ts ij}(x\com y)\,R_{\ts jk}(y\com z)=
R_{\ts jk}(y\com z)\,\Rt_{\ts ij}(x\com y)\,\Rt_{\ts ik}(x\com z)
\end{equation}
for any indices $i\com j\com k\in\{1\lcd\lt+l\ts\}$ such that
$i\le\lt<j\com k$ and $j\neq k\ts$.
By changing $x$ to $-x$ in (\ref{3.7}) and then
using (\ref{3.6}), we obtain the relation
\begin{equation}\label{3.8}
\ \,\qquad
\Rb_{\ts ij}(x\com y)\,\Rb_{\ts ik}(x\com z)\,R_{\ts jk}(y\com z)=
R_{\ts jk}(y\com z)\,\Rb_{\ts ik}(x\com z)\,\Rb_{\ts ij}(x\com y)\,.
\end{equation}
Similarly, when $i\com j\le\lt<k$ and $i\neq j$,
we obtain the relation
\begin{equation}\label{3.85}
\ \,\qquad
R_{\ts ji}(y\com x)\,\Rb_{\ts ik}(x\com z)\,\Rb_{\ts jk}(y\com z)=
\Rb_{\ts jk}(y\com z)\,\Rb_{\ts ik}(x\com z)\ts R_{\ts ji}(y\com x)\ts.
\end{equation}
To obtain the relation (\ref{3.85}), we have also used the identity
\begin{equation}\label{3.75}
R_{ij}(-x\com\ns-y)=R_{ji}(y\com x)\,.
\end{equation}
The relations (\ref{3.5})\com\ts(\ref{3.55}) and
(\ref{3.6}) to (\ref{3.85}) are equalities
of functions that take their values in the operator algebra $\Bll\,$.

Now take $\lt+l$ complex variables
$x_1\lcd x_{\ts\lt+l}\,$.
Consider the rational function of these variables, 
\begin{equation}\label{3.9}
\prod_{1\le i<j\le\lt}^{\longrightarrow}\,
R_{\ts ji}(\ts x_j\com x_i\ts)
\ts\ \cdot\!\!\!
\prod_{\substack{1\le i\le\lt\\\lt<j\le\lt+l}}^{\longrightarrow}\,\ts
\Rb_{\ts ij}(\ts x_i\com x_j\ts)
\ts\ \cdot\!\!\!\!\!\!
\prod_{\lt<i<j\le\lt+l}^{\longrightarrow}\!\ns
R_{\ts ij}(\ts x_i\com x_j\ts)
\end{equation}
taking values in the algebra $\Bll\,$. In each of the three
ordered products in (\ref{3.9}), 
the pairs $(i\com j)$ are ordered lexicographically. 
By using the relations (\ref{3.85}), the entire product (\ref{3.9})
can be rewritten as
$$
\prod_{1\le i\le\lt}^{\longleftarrow}\ \biggl(\ 
\prod_{\lt<j\le\lt+l}^{\longrightarrow}\,
\Rb_{\ts ij}(\ts x_i\com\ts x_j\ts)\biggr)
$$
\vglue-6pt
\begin{equation}\label{3.999}
\times\ 
\prod_{1\le i<j\le\lt}^{\longrightarrow}\,
R_{\ts ji}(\ts x_j\com x_i\ts)
\ts\ \cdot\!\!\!\!\!\!
\prod_{\lt<i<j\le\lt+l}^{\longrightarrow}\!\ns
R_{\ts ij}(\ts x_i\com x_j\ts)\ .
\end{equation}
By using the relations (\ref{3.8}), the product (\ref{3.9})
can also be rewritten as
$$
\prod_{1\le i<j\le\lt}^{\longrightarrow}\,
R_{\ts ji}(\ts x_j\com x_i\ts)
\ts\ \cdot\!\!\!\!\!\!
\prod_{\lt<i<j\le\lt+l}^{\longrightarrow}\!\ns
R_{\ts ij}(\ts x_i\com x_j\ts)
$$
\vglue-6pt
\begin{equation}\label{3.9999}
\times\ 
\prod_{1\le i\le\lt}^{\longrightarrow}\ \biggl(\ 
\prod_{\lt<j\le\lt+l}^{\longleftarrow}\,
\Rb_{\ts ij}(\ts x_i\com x_j\ts)\biggr)\ .
\end{equation} 
In the product displayed in the first line of
(\ref{3.999}), the factors corresponding to
$i=1\lcd\lt$ are arranged from right to left, as indicated
by the reversed arrow; a similar convention has been used in the
display (\ref{3.9999}).


\smallskip\medskip\noindent\textbf{3.3.}
Denote by $W_l$ and $W_{\lt\ts}^{\ts\ast}$ the tensor products of 
respectively
$l$ and $\lt$ copies of the vector spaces $\CC^L$ and $(\CC^L)^\ast$.
Choose two standard tableaux $\La$ and $\Lat\ts$, of shapes $\la$ and $\lat$
respectively. Denote by $F_\La$ and $F_{\Lat}$ the linear operators on 
respectively $W_l$ and $W_{\lt\ts}^{\ts\ast}$, 
corresponding to the diagonal matrix elements
$f_\La\in\CC S_l$ and $f_{\Lat}\in\CC S_{\ts\lt}\,\ts$; see Subsection 2.1. 
Let $P_{\ts\lt}$ be the linear operator on $W_{\lt\ts}^{\ts\ast}\ts$, 
reversing the order of $\lt$ tensor factors $(\CC^L)^\ast$.

Take $l$ complex variables $t_1(\La)\lcd t_l(\La)$ 
constrained by (\ref{1.22}). Setting
\begin{equation}\label{se}
x_{\ts\lt+k}=c_{k}(\La)+t_{k}(\La)
\text{\ \ for each\ \ $k=1\lcd l$}
\end{equation}
in the product
$$
\prod_{\lt<i<j\le\lt+l}^{\longrightarrow}\!\ns
R_{\ts ij}(\ts x_i\com x_j\ts)
$$
we obtain a rational function of the constrained variables
$t_1(\La)\lcd t_l(\La)\ts$. By Proposition 2.2, this function is regular at 
\begin{equation}\label{or}
t_1(\La)=\ldots=t_l(\La)=0\ts.
\end{equation}
The value of this function at (\ref{or})
coincides with the operator $1\ot F_\La$ on
\begin{equation}\label{del}
\Wll\,=\,W_{\lt\ts}^{\ts\ast}\ot W_l\,.
\end{equation}

Furthermore, take the $\lt$ constrained variables
$t_1(\Lat)\lcd t_{\ts\lt\ts}(\Lat)\ts$. Setting
\begin{equation}\label{set}
x_{\ts\lt-k+1}=c_{k}(\Lat)+t_{k}(\Lat)
\text{\ \ for each\ \ $k=1\lcd\lt$}
\end{equation}
in
$$
\prod_{1\le i<j\le\lt}^{\longrightarrow}\,
R_{\ts ji}(\ts x_j\com x_i\ts)
\ =\ 
P_{\ts\lt}\hskip4pt\cdot\hskip-6pt
\prod_{1\le i<j\le\lt}^{\longrightarrow}\,
R_{\ts ij}(\ts x_{\ts\lt-i+1}\com x_{\ts\lt-j+1}\ts)\hskip4pt\cdot\hskip3pt
P_{\ts\lt}
$$
we obtain a rational function of $t_1(\Lat)\lcd t_{\ts\lt\ts}(\Lat)$
which is regular at
\begin{equation}\label{ort}
t_1(\Lat)=\ldots=t_{\ts\lt\ts}(\Lat)=0\ts.
\end{equation}
Here we apply Proposition 2.2 to the tableau $\Lat$ instead of $\La\ts$,
and use the Yang-Baxter relation (\ref{3.5}) repeatedly.
The value of this restriction at (\ref{ort}) equals the operator
$(\ts P_{\ts\lt}\,F_{\Lat}\,P_{\ts\lt}\ts)\ot1$ 
on the vector space $\Wll\,$.

Finally, note that for any $k\in\{1\lcd l\}$ and $i\in\{1\lcd\lt\ts\}$
we have
$$
c_k(\La)+c_i(\Lat)\,\ge\,2-\lap_1-\lat^{\ts\prime}_1\,\ge\,2-L
$$
due to our assumptions on the partitions $\la$ and $\lat\ts$.
Using observations made in this and the previous subsection,
we arrive at the following proposition.

\begin{Proposition}
Restriction of the rational function 
{\rm(\ref{3.9})} to {\rm(\ref{se}),(\ref{set})} 
is regular at {\rm(\ref{or}),(\ref{ort})}.
The value of this restriction at {\rm(\ref{or}),(\ref{ort})} equals
$$
\prod_{1\le i\le\lt}^{\longrightarrow}\ \biggl(\ 
\prod_{1\le k\le l}^{\longrightarrow}\ 
\biggl(1-\frac{Q_{\ts\lt-i+1\ts,\ts\lt+k}}
{\ts c_k(\La)+{c_i(\Lat)}^{\phantom{\prime}}\ns\!+L}
\,\biggr)\ns\biggr)\ \cdot\ 
(\ts P_{\ts\lt}\,F_{\Lat}\,P_{\ts\lt}\ts)\ot F_\La
$$
$$
\ =\ (\ts P_{\ts\lt}\,F_{\Lat}\,P_{\ts\lt}\ts)\ot F_\La\ \,\cdot\,
\prod_{1\le i\le\lt}^{\longleftarrow}\ \biggl(\ 
\prod_{1\le k\le l}^{\longleftarrow}\ 
\biggl(1-\frac{Q_{\ts\lt-i+1\ts,\ts\lt+k}}
{\ts c_k(\La)+{c_i(\Lat)}^{\phantom{\prime}}\ns\!+L}
\,\biggr)\ns\biggr)\,.
$$
\end{Proposition}

\noindent
Let us define the linear operator $\FLa\ts$ on the vector space $\Wll$ 
as either of the two (equal) products displayed in Proposition 3.3. 
Using this
definition of $\FLa$ along with the relations (\ref{2.4}) and (\ref{3.8}),
we obtain the identity
$$
P_{\ts\lt+k\ts,\ts\lt+k+1}\ts 
R_{\ts\lt+k+1\ts,\ts\lt+k}\ts(\ts c_{k+1}(\La)\com c_k(\La))
\,\cdot\,\FLa
$$
\begin{equation}\label{kid}
\hskip42pt
\,=\ts\,F_{s_k\La,\Lat\ts}\,\cdot\,
P_{\ts\lt+k\ts,\ts\lt+k+1}\ts 
R_{\ts\lt+k\ts,\ts\lt+k+1}\ts(\ts c_{k}(\La)\com c_{k+1}(\La))
\end{equation}
for any $k\in\{1\lcd l-1\}$ such that the tableau $s_k\La$ is standard.
Similarly, by using the relation (\ref{3.85}) instead of (\ref{3.8}), 
we obtain the identity
$$
P_{\ts\lt-i\ts,\ts\lt-i+1}\ts 
R_{\ts\lt-i\ts,\ts\lt-i+1}\ts(\ts c_{i+1}(\Lat)\com c_i(\Lat))
\,\cdot\,\FLa
$$
\begin{equation}\label{iid}
\hskip48pt
\,=\ts\,F_{\La,s_i\Lat\ts}\,\cdot\,
P_{\ts\lt-i\ts,\ts\lt-i+1}\ts 
R_{\ts\lt-i+1\ts,\ts\lt-i}\ts(\ts c_{i}(\Lat)\com c_{i+1}(\Lat))
\end{equation}
for any index $i\in\{1\lcd\lt-1\}$ such that 
the tableau $s_i\Lat$ is standard.
The identities (\ref{kid}) and (\ref{iid}) provide relations
between the operators $\FLa$ for different pairs $\La$ and $\Lat$,
see the end of Subsection 2.1. 

Denote by $\beta$ the involutive antiautomorphism of the algebra $\Bll\,$,
such that each of the generators $P_{ij}$ and $Q_{ij}$ of $\Bll$  
is $\beta\ts$-invariant. The equality of two products displayed in
Proposition 3.3 can be then reformulated~as

\begin{Corollary}
The element\/ $\FLa\in\Bll\ts$ is\/ $\beta\ts$-invariant.
\end{Corollary}


\noindent\textbf{3.4.}
Denote by $\VLa$ the image of the operator
$\FLa$ on the vector space $\Wll\ts$.
Denote respectively by $V_\La$ and $V_{\ns\Lat}^{\ts\ast}$ the images of
the operators $F_\La$ and $F_{\Lat}$ on the vector spaces $W_l$ and 
$W_{\lt\ts}^{\ts\ast}\ts$. Then
\begin{equation}\label{inc}
\VLa\ts\subset\ts
(\ts P_{\ts\lt}\cdot V_{\Lat}^{\ts\ast})\ot V_\La
\end{equation}
due to the (second) definition of the operator $\FLa\,$. Take the subspace
$\Wllo\ts$ of traceless tensors. The next proposition is pivotal
for the present article.

\begin{Proposition}
We have the equality of vector spaces
$$
\VLa\ts=\ts
((\ts P_{\ts\lt}\cdot V_{\ns\Lat}^{\ts\ast})\ot V_\La)
\,\cap\,\Wllo\ts\,.
$$
\end{Proposition}

\begin{proof}
By definition, each of the operators $Q_{ij}$ vanishes on the subspace
$\Wllo\,$. So by the (second) definition of the operator $\FLa\ts$,
its action on~the subspace $\Wllo\ts$ coincides with that
of the operator $(\ts P_{\ts\lt}\,F_{\Lat}\,P_{\ts\lt}\ts)\ot F_\La\ts$.~Hence
$$
\VLa\,\supset
((\ts P_{\ts\lt}\,F_{\Lat}\,P_{\ts\lt}\ts)\ot F_\La\ts)\,\cdot\ts\Wllo
\ts\,=\,
((\ts P_{\ts\lt}\cdot V_{\ns\Lat}^{\ts\ast})\ot V_\La)
\,\cap\,\Wllo\,\ts.
$$
Due to (\ref{inc}), it now remains to show that $\VLa\subset\Wllo\,$.
Equivalently, we have to show that $Q_{ij}\,\FLa=0$
for any possible indices $i$ and $j\ts$.

By the (first) definition of $\FLa\,$, this operator
is divisible on the left by
$$
1-\frac{Q_{\ts\lt,\lt+1}}
{\ts c_1(\La)+{c_1(\Lat)}^{\phantom{\prime}}\ns\!+L}
\ =\ 1-\frac{\ts Q_{\ts\lt,\lt+1}\ts}{L}\,\ts.
$$
So the relation 
$$
Q_{\ts\lt,\lt+1}^{\ts2}\ns=\ts L\cdot Q_{\ts\lt,\lt+1}
$$ implies
the equality $Q_{ij}\,\FLa=0$ for $i=\lt\com j=\lt+1$ 
and arbitrary $\La\com\Lat\ts$.

If for some $k\in\{1\lcd l-1\}\ts$
the tableau $s_k\La$ is standard, then by (\ref{kid})
$$
P_{\ts\lt+k\ts,\ts\lt+k+1}\ts\FLa\,=\,
\frac{\FLa}{c_{k+1}(\La)-c_k(\La)}\ \ts+
$$
$$
F_{s_k\La,\Lat\ts}\cdot
\biggl(
P_{\ts\lt+k\ts,\ts\lt+k+1}-\frac1{\ts c_k(\La)-c_{k+1}(\La)\ts}
\ts\biggr)\ts.
$$
If the tableau $s_k\La$ is not standard, then 
$$
P_{\ts\lt+k\ts,\ts\lt+k+1}\ts\FLa\,=\,\pm\,\FLa\ts\,;
$$
the latter relation follows from the (second) definition of $\FLa$ 
by (\ref{2.3}).
If for some $i\in\{1\lcd l-1\}\ts$
the tableau $s_i\Lat$ is standard, then due to (\ref{iid})
$$
P_{\ts\lt-i\ts,\ts\lt-i+1}\ts\FLa\,=\,
\frac{\FLa}{{c_{i+1}(\Lat)-c_i(\Lat)}^{\phantom{\prime}}\ns\!}\ \ts+
$$
$$
F_{\La,s_i\Lat\ts}\cdot
\biggl(
P_{\ts\lt-i\ts,\ts\lt-i+1}-
\frac1{\ts{c_i(\Lat)-c_{i+1}(\Lat)}^{\phantom{\prime}}\ns\!\ts}
\ts\biggr)\,.
$$
If the tableau $s_i\Lat$ is not standard, then again due to (\ref{2.3})
$$
P_{\ts\lt-i\ts,\ts\lt-i+1}\ts\FLa\,=\,\pm\,\FLa\,.
$$

For any permutation $s\in\Sll\ts$, let $P_s$ be
the corresponding linear operator on the tensor product $\Wll\ts$.
The four relations displayed above imply that
for certain operators $R_{\La^\prime\ns\Lat^\prime}(s)$ on $\Wll\ts$, 
we have the equality
$$
P_s\,\FLa=\ts
\sum_{\La^\prime\ns,\Lat^\prime}\,\ts 
F_{\La^\prime\ns\Lat^\prime}\,R_{\La^\prime\ns\Lat^\prime}(s)\,;
$$
here $\La^\prime\com\Lat^{\ts\prime}$ range
over all standard tableaux of shapes $\la\ts\com\lat$ respectively.
For any $i\in\{1\lcd\lt\ts\}$ and $j\in\{\lt+1\lcd\lt+l\ts\}$
we can find a permutation $s\in\Sll$ so that $s(i)=\lt$ and
$s(j)=\lt+1$. Then we get the equalities
$$
Q_{ij}\,\FLa\,=\,P_s^{\ts-1}\ts Q_{\ts\lt,\lt+1}\,P_s\ts\FLa
\,=\,\sum_{\La^\prime\ns,\Lat^\prime}\,
P_s^{\ts-1}\ts Q_{\ts\lt,\lt+1}\ts F_{\La^\prime\ns\Lat^\prime}\ts
R_{\La^\prime\ns\Lat^\prime}(s)=0\quad\qed
$$
\end{proof}

\noindent
By definition, the associative algebra $\Bll$ is semisimple.
Let $\Bla$ be the simple ideal of $\Bll$
corresponding to the irreducible $\Bll\,$-module $U_{\lat}\ot U_\la\ts$.
Let $s_0$ be the element of
the maximal length in the symmetric group $S_{\lt}\ts$.
Proposition 3.4 provides the following characterization of $\FLa\ts$.

\begin{Corollary}
The operator $\FLa\in\Bll$
is the unique element of the simple ideal\/ $\Bla\,$,
with the image under\/ {\rm(\ref{3.22})} corresponding to the element
\begin{equation}\label{cor}
(\ts s_0\ts f_{\Lat}\ts s_0)\ts\cdot\ts\io_{\lt}\ts(f_\La)\,\in\,\CC\Sll\,.
\end{equation}
\end{Corollary}

\begin{proof}
By definition, we have 
$$
\FLa\ts\in\ts(\ts P_{\ts\lt}\,F_{\Lat}\,P_{\ts\lt}\ts)\ot F_\La+\,\Ill\,\ts.
$$
Hence the image of $\FLa$ under the 
homomorphism {\rm(\ref{3.22})} corresponds to 
the element (\ref{cor}).
Due to Proposition 3.4, we also have $\FLa\in\Bla\ts$.
But any element of the ideal $\Bla\subset\Bll$
is uniquely determined by the image of this element under
the homomorphism (\ref{3.22})\qed
\end{proof}

\noindent
Due to (\ref{2.01}), the above characterization of
the operator $\FLa$ implies that
$$
\FLa^{\,2}\,=\,\FLa\,\cdot\,
l\ts!\ts/\dim U_\la\,\cdot\,
\lt\ts!\ts/\dim U_{\lat}\,.
$$


\smallskip\noindent\textbf{3.5.}
We will now extend the results of Subsection 3.3 to standard tableaux
of skew shapes. Take any $m\in\{0\lcd l-1\}$ and $\mt\in\{0\lcd\lt-1\}$\ts. 
As in
Subsection~2.3, denote by $\Up$ be standard tableau obtained from $\La$
by removing the boxes with the numbers $m+1\lcd l$. Let $\mu$ be the shape of
the tableau $\Up$. Define a standard tableau $\Om$ of the skew shape
$\lm$ as in (\ref{Om}).
Similarly, define the partition $\mut$ of $\mt\ts$, the standard tableau $\Upt$
of shape $\mut\ts$, and the standard tableau $\Omt$ of shape $\lat/\mut\ts$.
Put
$$
n=l-m
\quad\text{and}\quad 
\nt=\lt-\mt\,.
$$

Take any $M\in\{0\lcd L-1\}$ and put $N=L-M$. 
Fix a decomposition $\CC^L=\CC^N\oplus\CC^M$,
it provides the decomposition $(\CC^L)^\ast=(\CC^N)^\ast\oplus(\CC^M)^\ast$.
Consider the tensor product
\begin{equation}\label{3.000}
\Wmm\ts=\ts(\CC^M)^\ast\ot\ldots\ot(\CC^M)^\ast\ot\CC^M\ot\ldots\ot\CC^M
\end{equation}
of $m$ copies of the defining representation $\CC^M$ of the group $GL_M\ts$,
and of $\mt$ copies of the representation $(\CC^M)^\ast\ts$.
We will use the factorization
$$
\Wmm\,=\,W_{\mt\ts}^{\ts\ast}\ot W_m\,.
$$
similar to (\ref{den}) and (\ref{del}).
Consider also the subspace of traceless tensors
$\Wmmo\subset\Wmm\,$.
There is a unique $GL_M\ts$-equivariant projector
$$
I_{m\mt}:\,\Wmm\,\to\,\Wmmo\,;
$$
see Subsection 3.1.
Assume that the irreducible rational representation $V_{\mu\ts\mut}$
of $GL_M$ appears~in $\Wmmo\subset\Wmm\,$. 
Equivalently, assume that $\mup_1+\mut_1^{\,\prime}\le M$.

In Section 1, we denoted respectively by $W_n$ and $W_{\ts\nt}^{\ts\ast}$ 
the tensor products of $n$ and $\nt$ copies of $\CC^N$ and $(\CC^N)^\ast$,
so that we had the factorization (\ref{den}). Let us now regard
the tensor product of vector spaces $W_{\ts\nt}^{\ts\ast}\ot\Wmm\ot W_n$ 
as a direct summand of $\Wll\ts$, 
using the decompositions $\CC^L=\CC^N\oplus\CC^M$
and $(\CC^L)^\ast=(\CC^N)^\ast\oplus(\CC^M)^\ast$ 
of the $\lt+l$ tensor factors in (\ref{3.00}). Let
$J_{m\mt}$ the projector of $\Wll$ to this direct summand.
Further, consider the subspace
\begin{equation}\label{3.41}
\Zmm\,=\,
W_{\ts\nt}^{\ts\ast}\ot\Wmmo\ot W_n
\,\subset\,
W_{\ts\nt}^{\ts\ast}\ot\Wmm\ot W_n
\,\subset\,
\Wll\,.
\end{equation}
Denote by $H_{m\mt}$ the composition of the operators $1\ot I_{m\mt}\ot1$
and $J_{m\mt}\ts$, 
this composition is a projector to the subspace $\Zmm\subset\Wll\,$.
For any linear operator $A$ on $\Wll\ts$
determine the operator on the same space,
\begin{equation}\label{ave}
A^\vee=\,H_{m\mt}\,A\,H_{m\mt}\,.
\end{equation}

We may also regard $A^\vee$ as an operator on the subspace $\Zmm\ts$.
For any $A\in\Bll$ let $\Th_{m\mt}(A)$ be the operator 
on the vector space $\Wmmo\ot W_{\ts\nt}^{\ts\ast}\ot W_n$,
obtained from $A^\vee$ by identifying this vector space with $\Zmm$
in~the natural way, that is by
exchanging the tensor factors $\Wmmo$ and $W_{\ts\nt}^{\ts\ast}$.
Let $\Bmm$ and $\Bnn$ be the commutants of the actions of the~groups
$GL_M$ and $GL_N$ on the vector spaces $\Wmm$ and $\Wnn$ respectively.
Consider also the two-sided ideal $\Imm\subset\Bmm$ 
defined as in Subsection 3.1.
The projector $H_{m\mt}$ commutes with action of subgroup
$GL_N\times GL_M\ts\subset\ts GL_L$ on $\Wll\,$. 
Hence for $A\in\Bll\ts$ 
$$
\Th_{m\mt}(A)\in(\,\Bmm/\,\Imm)\ot\Bnn\,;
$$
Here the quotient algebra $\ts\Bmm/\,\Imm\ts$ is identified with 
the image of $\CC S_{\mt m}$ in the operator algebra $\End(\Wmmo\ts)\ts$. 
Thus we obtain a linear map
$$
\Th_{m\mt}\ts:\,\Bll\,\rightarrow\,(\Bmm/\,\Imm)\ot\Bnn\,.
$$

As the mapping $A\mapsto A^\vee$ is linear, for any possible 
$i\com j\in\{1\lcd\lt+l\}$
$$
\Rp_{ij}(x\com y)=1-\frac{\Pp_{ij}}{x-y}
\qquad\textrm{and}\qquad
\Rbp_{ij}(x\com y)=1-\frac{\Qp_{ij}}{x+y+L}\,.
$$
Note that the operator $\Pp_{ij}$ on the subspace $\Zmm$ coincides
with restriction of the permutational operator $P_{ij}$ on $\Wll$
to this subspace, whenever
$$
j<i\le\nt
\quad\text{or}\quad
\nt<j<i\le\lt
\quad\text{or}\quad
\lt<i<j\le\lt+m
\quad\text{or}\quad
\lt+m<i<j\,.
$$
Also note that if
$$
i\le\nt
\quad\text{and}\quad
\lt+m<j\,,
$$ 
then the operator $\Qp_{ij}$ on the subspace
$\Zmm$ acts as (\ref{1.45}) in $i\ts$th tensor factor $(\CC^N)^\ast$
and the $(j\ns-\lt-m\ts)\ts$th tensor factor $\CC^N\ns$, and acts
as the identity in~all the remaining tensor factors.

Now consider the rational function (\ref{3.9}) with values
in the algebra $\Bll\ts$, of the variables $x_1\lcd x_{\ts\lt+l}\,$.
In the following proposition, these variables need not to be restricted
by (\ref{se}) and (\ref{set}).

\begin{Proposition}
The image of the product\/
{\rm(\ref{3.9})} under the map\/ $A\mapsto A^\vee$~is
$$
\prod_{1\le i<j\le\nt}^{\longrightarrow}
\Rp_{\ts ji}(\ts x_j\com x_i\ts)
\ts\ \cdot\hskip-12pt
\prod_{\substack{1\le i\le\nt \\ \lt+m<j\le\lt+l}}^{\longrightarrow}\!\!
\Rbp_{\ts ij}(\ts x_i\com x_j\ts)
\ts\ \cdot\hskip-16pt
\prod_{\lt+m<i<j\le\lt+l}^{\longrightarrow}\hskip-8pt
\Rp_{\ts ij}(\ts x_i\com x_j\ts)
$$
\begin{equation}\label{3.43}
\times
\prod_{\nt<i<j\le\lt}^{\longrightarrow}\,
\Rp_{\ts ji}(\ts x_j\com x_i\ts)
\ts\ \cdot\hskip-12pt
\prod_{\lt<i<j\le\lt+m}^{\longrightarrow}\!\!\!
\Rp_{\ts ij}(\ts x_i\com x_j\ts)\,.
\end{equation}
\end{Proposition}

\begin{proof}
Denote by $W_{\ts\lt+l}$ the tensor product of $\lt+l$ copies of
the vector space $\CC^L$. By conjugating
any linear operator $B$ on the vector space $W_{\ts\lt+l}$ relative to each
of the first $\lt$ tensor factors $\CC^L$, we obtain a certain linear
operator on $\Wll\,$, which will be denoted by $\Bt$.

As observed in Subsection 3.2, the ordered product (\ref{3.9}) can also
be written as (\ref{3.999}). The function (\ref{3.999}) 
of the variables $x_1\lcd x_{\ts\lt+l}$ takes its values
in the algebra of linear operators on $\Wll\,$.
This function can be obtained by applying the conjugation
$B\mapsto\Bt$ to the values of the function
$$
\prod_{1\le i<j\le\lt}^{\longleftarrow}\,
R_{\ts ji}(\ts x_j\com x_i\ts)
\ts\ \cdot\ns
\prod_{1\le i\le\lt}^{\longleftarrow}\ \biggl(\ 
\prod_{\lt<j\le\lt+l}^{\longleftarrow}\,
R_{\ts ij}(\ts x_i\com-\ts x_j\ns-\!L\ts)\biggr)
$$
\begin{equation}\label{3.99}
\times\!
\prod_{\lt<i<j\le\lt+l}^{\longrightarrow}\!\ns
R_{\ts ij}(\ts x_i\com x_j\ts)\ .
\end{equation}
In the entire product (\ref{3.99}), every factor is defined by
(\ref{3.45}) where $P_{ij}$ is the linear operator on $W_{\ts\lt+l}$
exchanging the $i\ts$th and $j\ts$th copies of $\CC^L$.
Using the Yang-Baxter relation (\ref{3.5}) repeatedly, the
product displayed in the second line of (\ref{3.99})
can be rewritten as 
$$
\prod_{\lt<i<j\le\lt+l}^{\longleftarrow}\!\ns
R_{\ts ji}(\ts-\ts x_j\ns-\!L\com-\ts x_i\ns-\!L\ts)\ .
$$

For the current proof only, let us reorder the
sequence of indices (\ref{ind}) as 
$$
\lt\com\lt-\ns1\lcd1\com\lt+l\lcd\lt+2\com\lt+1\,.
$$
We will employ the symbol $\prec$ when referring to this reordering.
Then order lexicographically the set of pairs of indices $(i\com j)$ 
with $i\prec j$ from the latter sequence; the lexicographical
ordering is taken with respect to $\prec$ here. The entire product
(\ref{3.99}) can then be rewritten as
\begin{equation}\label{product}
\prod_{\lt\ts\preccurlyeq\ts i\ts\prec\ts 
j\ts\preccurlyeq\ts\lt+1}^{\longrightarrow}\hskip-6pt
R_{\ts ij}(\ts y_i\com y_j\ts)
\end{equation}
where
$$
y_i\ts=\,
\left\{
\begin{array}{ll}
x_i
&\ \ \textrm{if}\ \ \ i=1\lcd\lt\,;\\[2pt]
-x_i\!-\!L
&\ \ \textrm{if}\ \ \ i=\lt+1\lcd\lt+l\,.
\end{array}
\right.
$$

The ordered product (\ref{product}) can be expanded as a linear combination
of the linear operators $P_s$ on $W_{\lt+l}$ corresponding
to permutations $s\in S_{\lt+l}\,$, with the coefficients
from the field $\CC(\ts y_1\lcd y_{\ts \lt+l}\ts)\ts$.
By the definition (\ref{ave}), we have 
$\Pt_s{\hskip-4.5pt}^\vee=0$ unless the permutation $s$ preserves each of
the two subsets
$$
\{\nt+1\lcd\lt\}
\ts\com\ts
\{\lt+1\lcd\lt+m\}
\subset
\{1\lcd\lt+l\}\,.
$$ 
Using Proposition 2.3 twice, we now establish that the image of
(\ref{product}) under the mapping $B\mapsto\Bt{}^\vee$
coincides with the image of
$$
\prod_{\lt\preccurlyeq i\ts\prec j\prec\nt}^{\longrightarrow}
R_{\ts ij}(\ts y_i\com y_j\ts)
\ \ts\cdot\hskip-8pt
\prod_{\nt\preccurlyeq i\prec j\prec\lt+m}^{\longrightarrow}\hskip-6pt
R_{\ts ij}(\ts y_i\com y_j\ts)
\ \ts\cdot\hskip-12pt
\prod_{\lt+m\preccurlyeq i\prec j\preccurlyeq\lt+1}^{\longrightarrow}
\hskip-8pt
R_{\ts ij}(\ts y_i\com y_j\ts)\,.
$$
The latter image is easy to describe, because the three ordered products~over
$$
\lt\preccurlyeq i\ts\prec j\prec\nt\,,\ \  
\nt\preccurlyeq i\prec j\prec\lt+m\ \ \text{and}\ \   
\lt+m\preccurlyeq i\prec j\preccurlyeq\lt+1
$$
in the above display pairwise commute, and the image of
the entire displayed product is multiplicative relative to
the factorization into these three. Thus the image of
(\ref{product}) under $B\mapsto\Bt{}^\vee$ equals
$$
\prod_{\nt\preccurlyeq i\preccurlyeq1}^{\longrightarrow}\ \biggl(\ 
\prod_{\lt+l\preccurlyeq j\prec\lt+m}^{\longleftarrow}\,
\Rbp_{\ts ij}(\ts x_i\com x_j\ts )\biggr)
\ \cdot\hskip-2pt
\prod_{\nt\preccurlyeq i\prec j\preccurlyeq1}^{\longleftarrow}
\Rp_{\ts ji}(\ts x_i\com x_j\ts)
\,\ \cdot\hskip-12pt
\prod_{\lt+l\preccurlyeq i\prec j\prec\lt+m}^{\longrightarrow}\hskip-8pt
\Rp_{\ts ji}(\ts x_j\com x_i\ts)
$$
$$
\times
\prod_{\lt\preccurlyeq i\ts\prec j\prec\nt}^{\longleftarrow}
\Rp_{\ts ij}(\ts x_i\com x_j\ts)
\ \ts\cdot\hskip-8pt
\prod_{\lt+m\preccurlyeq i\prec j\preccurlyeq\lt+1}^{\longrightarrow}
\hskip-8pt
\Rp_{\ts ji}(\ts x_j\com x_i\ts)\ \,=
$$
$$
\prod_{1\le i<j\le\nt}^{\longrightarrow}
\Rp_{\ts ji}(\ts x_j\com x_i\ts)
\ \ts\cdot\,\,
\prod_{1\le i\le\nt}^{\longrightarrow}\ \biggl(\ 
\prod_{\lt+m<j\le\lt+l}^{\longrightarrow}\hskip-3pt
\Rbp_{\ts ij}(\ts x_i\com x_j\ts )\biggr)
\ \ts\cdot\hskip-11pt
\prod_{\lt+m<i<j\le\lt+l}^{\longrightarrow}\hskip-9pt
\Rp_{\ts ij}(\ts x_i\com x_j\ts)
$$
$$
\times
\prod_{\nt<i<j\le\lt}^{\longrightarrow}
\Rp_{\ts ji}(\ts x_j\com x_i\ts)
\,\ \cdot\hskip-7pt
\prod_{\lt<i<j\le\lt+m}^{\longrightarrow}\hskip-5pt
\Rp_{\ts ij}(\ts x_i\com x_j\ts)\,\ts;
$$
here we used the analogues of the relations (\ref{3.5}) and (\ref{3.85}),
for the operators $\Rp_{\ts ij}(\ts x\com y\ts)$ and
$\Rbp_{\ts ij}(\ts x\com y\ts)$ instead of
of $R_{\ts ij}(\ts x\com y\ts)$ and $\Rb_{\ts ij}(\ts x\com y\ts)$.
At the right hand side of these equalities, we already employed
the natural ordering of indices (\ref{ind}), and the usual
lexicographical ordering of pairs of indices~\qed
\end{proof}

\noindent
Let us now make the substitutions (\ref{se}),(\ref{set}) and then the
substitutions (\ref{or}),(\ref{ort}) in Proposition 3.5. Note that 
$$
c_k(\La)\ts=\,
\left\{
\begin{array}{ll}
c_k(\ts\Up)
&\ \ \textrm{if}\ \ \ k\le m\ts,\\[2pt]
c_{k-m}(\Om)
&\ \ \textrm{if}\ \ \ k>m
\end{array}
\right.
$$
and
$$
c_k(\Lat)\ts=\,
\left\{
\begin{array}{ll}
c_k(\ts\Upt)
&\ \ \textrm{if}\ \ \ k\le\mt\ts,\\[2pt]
c_{k-\mt}(\Omt)
&\ \ \textrm{if}\ \ \ k>\mt\,.
\end{array}
\right.
$$
After the substitutions, the rational function of
$x_{\ts\nt+1}\lcd x_{\ts\lt+m}$
displayed in the second line line of (\ref{3.43}) takes the operator value
$$
1\ot((\,P_{\ts\mt}\,F_{\ts\Upt}\,P_{\ts\mt})\ot F_{\ts\Up})\,|\,\Wmmo\ts)
\ot1\,,
$$
where $F_{\ts\Up}$ and $F_{\ts\Upt}$ are operators on respectively
$W_m$ and $W_{\ts\mt}^{\ts\ast}$, while
$P_{\ts\mt}$ is the operator on 
$W_{\mt\ts}^{\ts\ast}$ reversing the order
of tensor factors. Here we applied 
Proposition 2.2 to the (non-skew) standard tableaux $\Up$
and $\Upt$ instead of $\La\ts$.

Now consider the operator $\FOm$ on $\Wnn$ as defined in Subsection 1.3.
The product displayed in the first line of (\ref{3.43}) can be rewritten as
$$
\prod_{1\le i\le\nt}^{\longleftarrow}\ \biggl(\ 
\prod_{\lt+m<j\le\lt+l}^{\longrightarrow}\,
\Rbp_{\ts ij}(\ts x_i\com\ts x_j\ts)\biggr)
$$
$$
\times\ 
\prod_{1\le i<j\le\nt}^{\longrightarrow}\,
\Rp_{\ts ji}(\ts x_j\com x_i\ts)
\ts\ \cdot\hskip-16pt
\prod_{\lt+m<i<j\le\lt+l}^{\longrightarrow}\!\ns
\Rp_{\ts ij}(\ts x_i\com x_j\ts)\ ;
$$
see Subsection 3.2. By using Proposition 3.3 and the definition of
$\Th_{m\mt}\ts$, after
the above substitutions in Proposition 3.5, we obtain the equality
\begin{equation}\label{3.4444}
\Th_{m\mt}(\FLa)=
\bigl(
(\ts P_{\ts\mt}\,F_{\ts\Upt}\,P_{\ts\mt})\ot F_{\ts\Up})\,|\,\Wmmo\ts
\bigr)
\ot\FOm\,.
\end{equation}
The linear operator $\FLa$ on $\Wll$ has been defined as
either of the two (equal) products displayed in Proposition 3.3.
When $m+\mt>0$, the equality (\ref{3.4444}) can be
regarded as an alternative definition of the operator $\FOm$ on $\Wnn\ts$.


\bigskip\noindent\textbf{3.6.}
In this subsection, we need to indicate the basic vector 
spaces $\CC^N\com\CC^L$ and $\CC^M$ corresponding to the tensor products
(\ref{mix}),(\ref{3.00}) and (\ref{3.000}). For this reason,
we will employ respectively
the symbols $\Wnn\ts(N)\com\Wll\ts(L)$ and $\Wmm\ts(M)$ 
to denote these tensor products. In the same style, we will write
$W_n(N)\com W_{\nt}^{\ts\ast}(N)$ and $\Wnno\ts(N)\ts$.
The linear operator (\ref{1.45}) on $(\CC^N)^\ast\ot\CC^N$
will be denoted by $Q(N)\ts$. Following this style of notation
throughout the article would make our formulas rather cumbersome,
so we confine this style to the present subsection.  
Recall that here we have $M=L-N$.

The image of the operator $\FOm$ on the vector space
$\Wnn\ts(N)$ has been denoted by $\VOm\,$. 
If $M=0$ then $\mu=\mut=(0\com0\ts,\ts\ldots\ts)\ts$ and
$\Om=\La\com\Omt=\Lat\ts$.
In this special case, we have the equality (\ref{1.4444}) by Proposition 3.3.
The next proposition holds for arbitrary $M\ge0\ts$.

\begin{Proposition}
If\/ $\Vlm\neq\{0\}$, then $\VOm\neq\{0\}$.
\end{Proposition}

\begin{proof}
Let us realize the irreducible
representation $\Vla$ of the group $GL_L$
from (\ref{1.0}) as the image $\VLa\subset\Wll\ts(L)$
of the operator $\FLa\,$. Choose a basis $v_1\lcd v_L$ in $\CC^L$ such
that the subspaces $\CC^N$ and $\CC^M$ of $\CC^L$ are spanned
respectively by $v_1\lcd v_N$ and $v_{N+1}\lcd v_L\,$.
Let $v_1^\ast\lcd v_L^\ast$ be the dual basis in $(\CC^L)^\ast$.
For any non-negative integers $k,\kt$ and $t$ such that $k+t\le l$ and
$\kt+t\le\lt\ts$, denote by $Z_{k\kt}^{\ts(t)}$ the subspace of 
$\Wll\ts(L)$ spanned by the vectors
$$
\sum_{N<a_1,\ldots,\ts a_t\le L}
u\ot 
v_{a_t}^{\ts\ast}\ot\ldots\ot v_{a_1}^{\ts\ast}\ot
w\ot
v_{a_1}\ot\ldots\ot v_{a_t}\ot 
v\,,
$$
for all $u\in W_{\ts\lt-\kt-t\ts}^{\ts\ast}(N)\ts$,
all $v\in W_{\ts l-k-t\ts}(N)\ts$, 
and all $w\in W_{\ts k\kt\ts}^{\ts\circ}(M)\ts$.

The vector space $\Wll\ts(L)$
is the sum of its subspaces, obtained from all the subspaces
of the form $Z_{k\kt}^{\ts(t)}$ 
by the permutational operators $P_s\ts$
for some $s\in\Sll\ts$. 
This sum may be not direct, cf.\ \cite[Section V.6]{W}.
However, every summand is preserved by 
the action of the subgroup $GL_M\subset GL_L\ts$.

For instance, consider the subspace $Z_{k\kt}^{\ts(t)}$ itself. 
Let $\bar{Z}_{k\kt}^{\ts(t)}$ be the subspace of $\Wll\ts(L)$ 
spanned, for the same $u,v$ and $w$ as above, by the vectors
$$
\sum_{1\le a_1,\ldots,\ts a_t\le N}
u\ot 
v_{a_t}^{\ts\ast}\ot\ldots\ot v_{a_1}^{\ts\ast}\ot
w\ot
v_{a_1}\ot\ldots\ot v_{a_t}\ot 
v\,.
$$
Then

\vskip-20pt
$$
\FLa\cdot Z_{k\kt}^{\ts(t)}\ts=\ts\FLa\cdot\bar{Z}_{k\kt}^{\ts(t)}\,\ts;
$$

\smallskip\noindent
this observation follows from from the equalities 
$\FLa\,Q_{ij}=0$ for the pairs
\begin{equation}\label{3.54}
(i\com j)=
(\ts\lt-\kt\com\lt+k+1)
\lcd
(\ts\lt-\kt-t+1\com\lt+k+t)\,;
\end{equation}
here $Q_{ij}$ are operators in $\Wll\ts(L)\ts$.
These operator equalities are implied by
Proposition 3.4 and Corollary~3.3. 
Note that $\bar{Z}_{k\kt}^{\ts(t)}\subset Z_{k\kt}^{\ts(0)}\ts$. 
Also note that
$$
{\rm Hom}_{\,GL_M}
\bigl(\,V_{\ts\mu\mut}\ts\com\FLa\cdot Z_{k\kt}^{\ts(0)}\,\bigr)
\neq\{0\}
\ \ \Rightarrow\ \ k=m\hskip.5pt\com\kt=\mt\,.
$$
The subspace $Z_{m\mt}^{\ts(0)}\subset\Wll\ts(L)$
coincides with the subspace $\Zmm\ts$, see (\ref{3.41}). 

By definition, the subspace $Z_{k\kt}^{\ts(t)}$ is contained in the 
direct summand
$$
W_{\ts\lt-\kt-t\ts}^{\ts\ast}(N)\ot
W_{\ts k+t\ts,\ts\kt+t\ts}(M)\ot
W_{\ts l-k-t\ts}(N)\,\subset\,\Wll\ts(L)\ts.
$$
Let $J_{k+t,\kt+t}$ be the projector in $\Wll\ts(L)$ to this direct summand.
Let $I_{k\kt}^{\ts(t)}$ be the operator on the tensor product
$W_{k+t\ts,\ts\kt+t}(M)\ts$, acting as
$Q(M)\ns/M$ on the $i\ts$th and $j\ts$th
tensor factors for every pair of indices
$$
(i\com j)=
(\ts t\com\kt+k+t+1)
\lcd
(1\com\kt+k+2t)\,,
$$
and acting as

\vskip-16pt
$$
I_{k\kt}:\,\Wkk\ts(M)\,\to\,\Wkko\ts(M)
\nopagebreak
$$
on the remaining $\kt+k$ tensor factors of $W_{k+t\ts,\ts\kt+t}(M)\ts$;
see Subsection 3.5.

Let $H_{k\kt}^{(t)}$ be
the composition of operators $1\ot I_{k\kt}^{(t)}\ot1$
and $J_{k+t,\kt+t}$ on $\Wll\ts(L)\ts$, 
this composition is a $GL_M$-equivariant
projector to the subspace~$Z_{k\kt}^{\ts(t)}$.
In a similar way, by using $Q(N)\ns/N$ instead of $Q(M)\ns/M$, define
a projector 
$$
\bar{H}_{k\kt}^{\ts(t)}:\,\Wll(L)\,\to\,\bar{Z}_{k\kt}^{\ts(t)}\,.
$$
Using the notation of Subsection 3.5, we have the equality
$H_{k\kt}^{(0)}=H_{m\mt}\ts$.

If an irreducible representation of
$GL_M$ equivalent to $V_{\ts\mu\mut}$ occurs in the projection
$H_{k\kt}^{(t)}\cdot\VLa\ts$, 
it also occurs in the projection 
$\bar{H}_{k\kt}^{(t)}\cdot\VLa\ts$.
This follows from the equalities 
$Q_{ij}\ts F_\La=0$ for the pairs of indices $i$ and $j$ given by 
(\ref{3.54}). Either occurence implies that $k=m$.
But again, we have
$$
\bar{H}_{m\mt}^{(t)}\cdot\VLa\,\subset\,
H_{m\mt}^{(0)}\cdot\VLa\,=\,
H_{m\mt}\cdot\VLa\,.
$$

Now suppose that $\Vlm\neq\{0\}$. 
By the above argument, then
there exist a pair of permutations $s^{\ts\prime}$ and $s^{\ts\prime\prime}$ 
in $\Sll$ such that
\begin{equation}\label{3.56}
{\rm Hom}_{\,GL_M}(\,V_{\ts\mu\mut}\ts\com\ts 
H_{m\mt}\,P_{s'}\,\FLa\,P_{s''}\cdot\Zmm\ts)\,\neq\,\{0\}\,.
\end{equation}

Let $\Lap,\Lapp$ and $\Latp,\Latpp$ range over the sets of
standard tableaux of shapes $\la$ and $\lat$ respectively. Let 

\vskip-20pt
\begin{equation}\label{FLL}
F_{\La''\ns\Lat''}^{\ts\La'\Lat'}\,\in\,\Bll
\end{equation}
be the unique element of the simple ideal
$\Bla\subset\Bll\,$ with the image under the homomorphism 
{\rm(\ref{3.22})} corresponding to the element 
$$
(\ts s_0\ts f_{\Lat''}^{\ts\Lat'}\ts s_0)
\ts\cdot\ts
\io_{\lt}\ts(\ts f_{\La''}^{\ts\La'})\,\in\,\CC\Sll\,;
$$
see (\ref{fll}). If $\Lap=\Lapp$ and $\Latp=\Latpp$, 
then we have the equality
$$
F_{\La''\ns\Lat''}^{\ts\La'\Lat'}\,=\,F_{\La'\ns\Lat'}
$$
due to Corollary 3.4. 
The product $P_{s'}\,\FLa\,P_{s''}$ in (\ref{3.56}) can be 
written as a linear combination of the operators (\ref{FLL}), 
with the coefficients from $\CC\ts$. By (\ref{3.56}), 
there exists standard tableaux $\Lap,\Lapp$ and $\Latp,\Latpp$ such that
\begin{equation}\label{3.57}
{\rm Hom}_{\,GL_M}(\,V_{\ts\mu\mut}\ts\com\ts 
H_{m\mt}\,F_{\La''\ns\Lat''}^{\ts\La'\Lat'}\cdot\Zmm\ts)\,\neq\,\{0\}\,.
\end{equation}

By definition, the operator (\ref{FLL}) is divisible by
$(\ts P_{\ts\lt}\,F_{\Lat'}\ts P_{\ts\lt}\ts)\ot F_{\La'}$ 
on the left, and by
$(\ts P_{\ts\lt}\,F_{\Lat''}\ts P_{\ts\lt}\ts)\ot F_{\La''}$
on the right.
Now consider the tableaux $\Up^{\ts\prime}$ and $\Up^{\ts\prime\prime}$,
obtained by removing the boxes with the numbers $m+1\lcd l$ from the 
tableaux $\Lap$ and $\Lapp$ respectively.
Also consider the
the tableaux $\Upt^{\ts\prime}$ and $\Upt^{\ts\prime\prime}$,
obtained by removing the boxes with the numbers $\mt+1\lcd\lt$ from the 
tableaux $\Latp$ and $\Latpp$. Then
the restriction of the operator
$$
H_{m\mt}\,F_{\La''\ns\Lat''}^{\ts\La'\Lat'}
$$
to the subspace $\Zmm\subset\Wll(L)$ is divisible on the left by 
$$
1\ot((\,P_{\ts\mt}\,F_{\ts\Upt'}\,P_{\ts\mt})\ot F_{\ts\Up'})
\,|\,\Wmmo(M))
\ot1\,,
$$
cf.\ the proof of Lemma 2.3. Similarly, it
is divisible on the right by 
$$
1\ot((\,P_{\ts\mt}\,F_{\ts\Upt''}\,P_{\ts\mt})\ot F_{\ts\Up''})
\,|\,\Wmmo(M))
\ot1\,.
$$
Now the inequality (\ref{3.57}) implies that the shapes of the tableaux
$\Up^{\ts\prime}$ and $\Up^{\ts\prime\prime}$ coincide.
The same inequality implies that the shapes of
$\Upt^{\ts\prime}$ and $\Upt^{\ts\prime\prime}$ coincide as well.
Then 
$$
F_{\La''\ns\Lat''}^{\ts\La'\Lat'}\ts=\,F_{\ts\La'\ns\Lat'}\,F
$$
for the operator $F$ on $\Zmm$ corresponding to
some invertible element $f$ of the group ring of
$$
S_{\ts\nt\mt}\times\io_{\lt}\ts(S_{mn})\,\subset\,S_{\lt+l}\,,
$$
see (\ref{2.3}). Then
$$
F_{\La''\ns\Lat''}^{\ts\La'\Lat'}\cdot\Zmm\,=\,F_{\ts\La'\ns\Lat'}\cdot\Zmm\,.
$$
By applying now the relation (\ref{3.4444}) to the tableaux $\Lap$ and
$\Latp$ instead of $\La\ts$ and $\Lat\ts$, the inequality (\ref{3.57})
implies that the tableaux $\Up^{\ts\prime}$ and $\Upt^{\ts\prime}$
are of shapes $\mu$ and $\mut$ respectively.
Moreover, the left hand side of (\ref{3.57}) then
equals $V_{\Om'\ns\Omt'}$ for some standard tableaux $\Om^{\ts\prime}$
and $\Omt^{\ts\prime}$ of skew shapes $\lm$ and $\lmt$ respectively.
The inequality $V_{\Om'\ns\Omt'}\neq\{0\}$ implies that $\VOm\neq\{0\}$\qed
\end{proof}

\noindent
The vector space $\Vlm$ as defined by (\ref{1.0}),
comes with a natural action of the subgroup $GL_N\subset GL_L\ts$.
On the other hand,
the subspace $\VOm\subset\Wnn(N)$ is preserved by the action of $GL_N\ts$, 
because $\FOm\in\Bnn\ts$. Let us consider $\Vlm$ and $\VOm$ as
representations of the group $GL_N$. In Subsection 5.4
we will prove that
these representations are equivalent, this is what Proposition~1.3 states.
We will use Proposition~3.6 along with
an explicit description of the action of the subalgebra
$\AMN\subset\UMN\ts$ on the vector space $\VOm\,$.


\section{\hskip-3.5pt. Yangian representations}


\noindent\textbf{4.1.}
In this section we will work with the space of mixed tensors (\ref{mix}).
Thus our basic vector space will be $\CC^N$ like it was in Section 1,
and unlike it was in Section 3.  
We will also use the tensor product of vector spaces
\begin{equation}\label{aux}
\CC^N\ns\ot\Wnn\ts=\ts
\CC^N\ns\ot(\CC^N)^\ast\ot\ldots\ot(\CC^N)^\ast\ot\CC^N\ns\ot\ldots\ot\CC^N
\end{equation}
where the first tensor factor $\CC^N$ shall not be regarded
as a representation of the group $GL_N$, it will play only an auxiliary role.
The $1\ns+\ts\nt+n\ts$ tensor factors in (\ref{aux}) will be
labeled by the indices
$$
0\com1\lcd\nt\com\nt+\ns1\lcd\nt+n
$$
respectively from the left to right. 
For two distinct indices $i$ and $j$ such that
$$
0<i\com j\le\nt
\quad\text{or}\quad
\nt<i\com j
\quad\text{or}\quad
i=0
\quad\text{and}\quad  
\nt<j
$$
let $P_{ij}$ be the linear operator on (\ref{aux}), exchanging
the $i\ts$th and $j\ts$th tensor factors. Then determine 
the function $R_{ij}(x,y)$ by the same equality (\ref{3.45}) as before.
But now $R_{ij}(x,y)$ takes values in the algebra of operators on
(\ref{aux}).

Further, for any two indices
$i$~and~$j$ such that $0<i\le\nt<j\ts$, let
$Q_{ij}$ the operator on (\ref{aux}) acting as (\ref{1.45})
on the $i\ts$th and $j\ts$th tensor factors,
and acting as the identity on remaining 
$\nt\ns+n\ns-\!1\ts$ factors. For such $i$ and $j$ put
\begin{equation}\label{instead}
\Rb_{ij}(x\com y)\,=\,1-\frac{Q_{ij}}{x+y+N}\,\ts;
\end{equation}
note the difference between the equality displayed above, 
and the second equality in (\ref{3.555}). 
Furthermore, for any index $i$ such that $0<i\le\nt\ts$, put
$$
\Ra_{\ts0j}(x,y)=
\biggl(1+\frac{Q_{0i}}{x-y}\ts\biggr)
\cdot
\biggl(1-\frac1{(x-y)^{\ts2}}\biggr)^{\!-1}
$$
where $Q_{0i}$ is the linear 
operator on (\ref{aux}) acting as (\ref{1.45})
on the $i\ts$th and $0\ts$th tensor factors,
and acting as the identity on remaining 
$\nt\ns+n\ns-\!1\ts$ factors.

Note that for all possible indices $i\com j>0$ we have the relations
\begin{equation}\label{4.001}
R_{\ts0i}(x\com y)\,R_{\ts0j}(x\com z)\,R_{\ts ij}(y\com z)=
R_{\ts ij}(y\com z)\,R_{\ts0j}(x\com z)\,R_{\ts0i}(x\com y)\,,
\end{equation}
\begin{equation}\label{4.002}
\Ra_{\ts0i}(x\com y)\,\Ra_{\ts0j}(x\com z)\,R_{\ts ij}(y\com z)=
R_{\ts ij}(y\com z)\,\Ra_{\ts0j}(x\com z)\,\Ra_{\ts0i}(x\com y)\,,
\end{equation}
\begin{equation}\label{4.003}
\Ra_{\ts0i}(x\com-y)\,R_{\ts0j}(x\com z)\,\Rb_{\ts ij}(y\com z)\,=\,
\Rb_{\ts ij}(y\com z)\,R_{\ts 0j}(x\com z)\,\Ra_{\ts 0i}(x\com-y)\,.
\end{equation}
Here (\ref{4.001}) is a Yang-Baxter relation, while
(\ref{4.002}) and (\ref{4.003}) are derived from Yang-Baxter
relations, by using (\ref{3.55}) and (\ref{3.6}) where $L$ is replaced by $N$. 

For any $y_1\lcd y_{\ts\nt}\in\CC$ and $z_1\lcd z_n\in\CC$ 
consider the tensor product 
\begin{equation}\label{mod}
V(y_{\ts\nt})^\ast\ns\ot\ldots\ot V(y_1)^\ast\ns\ot 
V(z_1)\ot\ldots\ot V(z_n)
\end{equation}
of evaluation and dual evaluation $\YN\ts$-modules, these have
been defined in Subsection~1.4. The underlying space of
the $\YN\ts$-module (\ref{mod}) is $\Wnn\,$. Denote by $\rho_{n\nt}$ 
the homomorphism $\YN\to\End(\Wnn)$ corresponding to (\ref{mod}).
Consider the element (\ref{1.71}). In what follows, we identify the algebras
\begin{equation}\label{twoalg}
\End(\CC^N)\ot\End(\Wnn)=
\End(\ts\CC^N\ns\ot\Wnn)\,.
\end{equation}

\begin{Proposition}
Under the homomorphism
$$
\id\ot\rho_{n\nt}:\,
\End(\CC^N)\ot\YN
\,\to\,
\End(\ts\CC^N\ns\ot\Wnn)\,,
$$
$$
T(x)\,\mapsto\,
\Ra_{\ts01}(x\com y_{\ts\nt})\ts\ldots\ts\Ra_{\ts0\nt}(x\com y_1)\,
R_{\ts0,\nt+1}(x\com z_1)\ts\ldots\ts R_{\ts0,\nt+n}(x\com z_n)\,.
$$
\end{Proposition}

\begin{proof}
Firstly, consider the case when $n=1$ and $\nt=0$. Using the
matrix units $E_{ab}\in\End(\CC^N)$, the element
$P_{\ts01}\in\End(\ts\CC^N\ns\ot\CC^N)$ expands as 
$$
P_{\ts01}\,=
\sum_{a,b=1}^N\, E_{ab}\ot E_{ba}\,. 
$$
Using definitions (\ref{1.52}) and (\ref{tau}), we prove that 
under the homomorphism 
$$
\id\ot\rho_{10}:\,
\End(\CC^N)\ot\YN
\,\to\,
\End(\ts\CC^N\ns\ot\CC^N)\,,
$$
we have $T(x)\mapsto R_{\ts01}(x\com z_1)$ as required.

Secondly, consider the case when $n=0$ and $\nt=1$.
The operator $Q_{01}$ on $\CC^N\ns\ot(\CC^N)^\ast$ is obtained
from the operator $P_{\ts01}$ on $\CC^N\ns\ot\CC^N$ by conjugation in the
second tensor factor. Using our result from the previous case and
the relation (\ref{3.55}), we prove that
$T(x)\mapsto\Ra_{\ts01}(x\com y_1)$ 
under the homomorphism 
$$
\id\ot\rho_{01}:\,
\End(\CC^N)\ot\YN
\,\to\,
\End(\ts\CC^N\ns\ot(\CC^N)^\ast)\,.
$$
For arbitrary $n$ and $\nt$, Proposition~4.1 follows from its
two particular cases considered here, and from the definition 
(\ref{1.33}) of coproduct on $\YN$
\qed
\end{proof}

\noindent
Consider also the product of evaluation and
dual evaluation $\YN\ts$-modules
\begin{equation}\label{modrev}
V(z_n)\ot\ldots\ot V(z_1)\ns\ot
V(y_1)^\ast\ns\ot\ldots\ot V(y_{\ts\nt})^\ast\,.
\end{equation}
Here the ordering of the tensor factors is opposite
to that used in (\ref{mod}).
The underlying vector space of the $\YN\ts$-module (\ref{modrev})
can be identified with that of the $\YN\ts$-module (\ref{mod}) 
by using the permutational operator $(\ref{Pnn})$. Let 
$\sigma_{n\nt}:\YN\to\End(\Wnn)$ be the composition of
the homomorphism corresponding to the $\YN\ts$-module (\ref{modrev}), with
the conjugation by $(\ref{Pnn})$.
Our proof of Proposition 4.1 then implies the following result.

\begin{Corollary}
Under the homomorphism
$$
\id\ot\si_{n\nt}:\,
\End(\CC^N)\ot\YN
\,\to\,
\End(\ts\CC^N\ns\ot\Wnn)\,,
$$
$$
T(x)\,\mapsto\,
R_{\ts0,\nt+n}(x\com z_n)\ts\ldots\ts R_{\ts0,\nt+1}(x\com z_1)\,
\Ra_{\ts0\nt}(x\com y_1)\ts\ldots\ts\Ra_{\ts01}(x\com y_{\ts\nt})\,.
$$
\end{Corollary}

Now consider the standard action of the enveloping algebra
$\UN$ on the space of mixed tensors $\Wnn\,$.
Let us denote by $\varpi_{\ts n\nt}$ 
the corresponding homomorphism $\UN\to\End(\Wnn)\,$.
Using definitions (\ref{1.71})~and~(\ref{1.52}),
we obtain the next result.

\begin{Lemma}
Under the homomorphism
$$
\id\ot(\ts\varpi_{\ts n\nt}\circ\ts\pi_N):\,\ts
\End(\CC^N)\ot\YN
\,\to\,
\End(\ts\CC^N\ns\ot\Wnn)\,,
$$
$$
T(x)\,\mapsto\,1+
\bigl(\ts 
Q_{01}+\ldots+Q_{0\nt}-P_{\ts0,\nt+1}-\ldots-P_{\ts0,\nt+n}
\ts\bigr)\,\ts
x^{-1}\,.
$$
\end{Lemma}


\smallskip\noindent\textbf{4.2.}
Take any two standard tableaux $\Om$ and $\Omt\ts$, with respectively
$n$ and $\nt$ boxes.
Consider the operator $\FOm$ on the vector space $\Wnn\,$,
as defined in Subsection 1.3. We denoted
by $P_{\ts\nt+n}$ the operator on $\Wnn\ts$ reversing the order of
all $\nt+n$ tensor factors, see (\ref{Pnn}).
Let us prove Proposition~1.4.~Set
\begin{equation}\label{setz}
z_k=c_k(\Om)
\text{\ \ for each\ \ }
k=1\lcd n\ts;
\end{equation}
\begin{equation}\label{sety}
y_k=-\ts c_k(\Omt)\ns-\!M
\text{\ \ for each\ \ }
k=1\lcd\nt\ts.
\end{equation}
In these settings, Proposition 4.1 describes 
the action of the algebra~$\YN$ in the tensor product
of evaluation and dual evaluation modules, displayed
in \text{Proposition~1.4} at the bottom.
The underlying space of the tensor product
of evaluation and dual evaluation modules displayed
in Proposition~1.4 at the top, can be identified with $\Wnn$
by using the operator (\ref{Pnn}). 
The action of $\YN$ in the resulting module is now described 
by Corollary~4.1.

Now take $\nt+n$ complex variables $x_1\lcd x_{\ts\nt+n}\,$.
Let $F(\ts x_1\lcd x_{\ts\nt+n})$ be the rational function of
these variables, defined as the product
$$
\prod_{1\le i<j\le\nt}^{\longrightarrow}
R_{\ts ji}(\ts x_j\com x_i\ts)
\ \,\cdot\hskip-8pt
\prod_{\substack{1\le i\le\nt\\\nt<j\le\nt+n}}^{\longrightarrow}
\Rb_{\ts ij}(\ts x_i\com x_j\ts)
\ts\ \cdot\hskip-13pt
\prod_{\nt<i<j\le\nt+n}^{\longrightarrow}\hskip-6pt
R_{\ts ij}(\ts x_i\com x_j\ts)\,.
$$
According to notation of Subsection 4.1, the values of this 
function
belong to the algebra (\ref{twoalg}). Using the relations (\ref{3.75}) and
(\ref{4.001}) to (\ref{4.003}) repeatedly, 
we obtain the equality of rational functions with values in 
the algebra (\ref{twoalg}),
$$
\Ra_{\ts01}(\ts x\com-x_1)
\ts\ldots\ts
\Ra_{\ts0\nt}(\ts x\com-x_{\ts\nt})\,
R_{\ts0,\nt+1}(\ts x\com x_{\ts\nt+1})
\ts\ldots\ts
R_{\ts0,\nt+n}(\ts x\com x_{\ts\nt+n})
$$
$$
\times\ 
F(\ts x_1\lcd x_{\ts\nt+n})
\,=\,
F(\ts x_1\lcd x_{\ts\nt+n})
\ \times
$$
$$
R_{\ts0,\nt+n}(\ts x\com x_{\ts\nt+n})
\ts\ldots\ts
R_{\ts0,\nt+1}(\ts x\com x_{\ts\nt+1})\,
\Ra_{\ts0\nt}(\ts x\com-x_{\ts\nt})
\ts\ldots\ts
\Ra_{\ts01}(\ts x\com-x_1)\,.
$$
In this equality, let us substitute
$$
x_{\ts\nt+k}=c_{k}(\Om)+t_{k}(\Om)
\text{\ \ for each\ \ $k=1\lcd n$}
$$
where the complex variables $t_1(\Om)\lcd t_n(\Om)$
are constrained by (\ref{1.0002}). By using the variables
$t_1(\Omt)\lcd t_{\ts\nt}(\Omt)$ constrained similarly, also substitute
$$
x_{\ts\nt-k+1}=c_{k}(\Omt)+M+t_{k}(\Omt)
\text{\ \ for each\ \ $k=1\lcd\nt\ts$.}
$$
We then obtain an equality of rational functions of the
constrained variables, and of the variable $x\ts$. 
By setting in the resulting equality 
$$
t_1(\Om)=\ldots=t_n(\Om)=0
\quad\text{and}\quad
t_1(\Omt)=\ldots=t_{\ts\nt}(\Omt)=0\,,
$$
we obtain an equality of rational functions
of the variable $x$ only. 
In view of Proposition 4.1 
and Corollary 4.1, 
the last equality proves Proposition~1.4.
Indeed, after all the above substitutions 
$F(\ts x_1\lcd x_{\ts\nt+n})$ takes the value
$1\ot\FOm\ts$ in the algebra (\ref{twoalg}), see 
the definition (\ref{maindef}).
Here we have used the equality between (\ref{3.9}) and (\ref{3.999}),
established in Subsection 3.2. Note that
the equality between (\ref{3.9}) and (\ref{3.9999}) now implies
Lemma 1.3.


\smallskip\medskip\noindent\textbf{4.3.}
Consider the embedding $\UN\to\YN$ as defined by ({\ref{4.4}).
By the definition (\ref{1.52}),
the homomorphism $\pi_N:\YN\to\UN$ is identical on the subalgebra 
$\UN\subset\YN\ts$.
The automorphism $\om_N$ of $\YN$
is also identical on this subalgebra, see the definition (\ref{1.51}).
It follows that the restriction of $\pi_{NM}$ to the subalgebra
$\UN\subset\YN$ coincides with the natural embedding $\UN\to\UMN\ts$.
For the particular choice (\ref{1.62}) of the series $g(x)$,
the automorphism (\ref{1.61}) is identical on 
$\UN\subset\YN\ts$, because the coefficient at 
$x^{-1}$ in the expansion of (\ref{1.62}) is $0\ts$.
So the action of the subalgebra
$\UN\subset\YN\ts$ in the $\YN\ts$-module $\Vlm$ coincides
with its natural action, corresponding to the natural
action of $GL_N$ in $\Vlm$.

Further, the action of the subalgebra
$\UN\subset\YN\ts$ in any evaluation module $V(z)$ over $\YN\ts$,
and in any dual evaluation module $V(z)^\ast$ over $\YN\ts$,
coincides with the natural action of $\UN$
in the underlying vector space $\CC^N$ of these modules.
The assigment (\ref{4.4}) determines a Hopf algebra
embedding, see (\ref{Hopfemb}).
Hence the action of the subalgebra
$\UN\subset\YN\ts$ in the tensor product of 
$\YN\ts$-modules displayed in Proposition 1.4 at the bottom,
coincides with the natural action of $\UN$
in the vector space $\Wnn\ts$. The $\YN\ts$-module $\VOm$ was
defined as a submodule of that tensor product.
Thus the action of $\UN\subset\YN\ts$ in this submodule
coincides with the natural action of $\UN$ in the subspace
$\VOm\subset\Wnn\ts$.


\smallskip\medskip\noindent\textbf{4.4.}
Take the $\YN\ts$-module (\ref{mod}) in the particular case
when $\nt=0$~and 
\begin{equation}\label{setext}
z_k=z-k+1
\text{\ \ for each\ \ }
k=1\lcd n\ts.
\end{equation}
Here $z$ is arbitrary complex papameter. The underlying vector
space of this module is the tensor product $W_n$ of $n$ copies of $\CC^N$.
Let $A_n$ be the operator of antisymmetrization on $W_n\ts,$ 
normalized so that $A_n^2=A_n\ts$. In the notation of
Subsection 4.1 we then have the equality in the algebra (\ref{twoalg}),
$$
\prod_{1\le i<j\le n}^{\longrightarrow}
R_{\ts ij}(\ts z_i\com z_j\ts)
\ =\ n\ts!\,\cdot\,1\ot A_n\,;
$$
see the end of Subsection 1.2.
The arguments of Subsection 4.2 now imply
$$
R_{\ts01}(\ts x\com z_1)
\ts\ldots\ts
R_{\ts0n}(\ts x\com z_n)
\,\cdot\ts(\ts1\ot A_n) 
$$
\begin{equation}\label{ext}
=\ 
(\ts1\ot A_n)\,\cdot\ts
R_{\ts0n}(\ts x\com z_n)
\ts\ldots\ts
R_{\ts01}(\ts x\com z_1)\,.
\end{equation}
Due to Propositition 4.1, the equality (\ref{ext}) indicates
that the $n\ts$th exterior power of $\CC^N$ is a submodule
of the $\YN\ts$-module (\ref{mod}), for $\nt=0$ and 
the parameters $z_1\lcd z_n$ defined by (\ref{setext}).
Denote by $V_n(z)$ this submodule.

Now consider the $\YN\ts$-module (\ref{mod}) in the case
when $n=0$ and 
\begin{equation}\label{setdext}
y_{\ts\nt-k+1}=z-k+1
\text{\ \ for each\ \ }
k=1\lcd\nt\ts.
\end{equation}
The underlying vector
space of this module is
$W_{\ts\nt}^{\ts\ast}\ts$.
Let $A_{\ts\nt}$ be the normalized operator of antisymmetrization
on $W_{\ts\nt}^{\ts\ast}\ts$. Similarly to (\ref{ext}), we have
$$
\Ra_{\ts01}(\ts x\com y_{\ts\nt})
\ts\ldots\ts
\Ra_{\ts0\hskip.5pt\nt}(\ts x\com y_1)
\,\cdot\ts(\ts1\ot A_{\ts\nt}) 
$$
\begin{equation}\label{dext}
=\ 
(\ts1\ot A_{\ts\nt})\,\cdot\ts
\Ra_{\ts0\hskip.5pt\nt}(\ts x\com y_1)
\ts\ldots\ts
\Ra_{\ts01}(\ts x\com y_{\ts\nt})\,.
\end{equation}
Due to Corollary 4.1, the equality (\ref{dext}) shows
that the $\nt\ts$th exterior power of $(\CC^N)^\ast$ is a submodule
of the $\YN\ts$-module (\ref{mod}), for $n=0$ and for
the parameters $y_1\lcd y_{\nt}$ defined by (\ref{setdext}).
Denote by $V_{\ts\nt}(z)^{\ts\prime}$ this submodule.

\begin{Lemma}
The\/ $\YN\ts$-module $V_n(z)^\ast$ dual to $V_n(z)\ts$, is equivalent to
$V_n(z)^{\ts\prime}$.
\end{Lemma}

\begin{proof}
The action of $\YN$ on $V_n(z)$ can be described by the assignment
$$
T(x)\,\mapsto\,
R_{\ts01}(\ts x\com z_1)
\ts\ldots\ts
R_{\ts0n}(\ts x\com z_n)
\,\cdot\ts(\ts1\ot A_n)\,,
$$ 
cf.\ Proposition 4.1.
The same action can also be described by the assignment
$$
T(x)^{-1}\,\mapsto\,
R_{\ts0n}(\ts x\com z_n)^{-1}
\ts\ldots\ts
R_{\ts01}(\ts x\com z_1)^{-1}
\,\cdot\ts
(\ts1\ot A_n)
$$
$$
=\ 
(\ts1\ot A_n)
\,\cdot\ts
R_{\ts01}(\ts x\com z_1)^{-1}
\ts\ldots\ts
R_{\ts0n}(\ts x\com z_n)^{-1}\,,
$$
see (\ref{ext}).
Using the definition (\ref{antip}) of the antipode ${\rm S}$
along with (\ref{3.55}), the action of $\YN$ 
on the module dual to $V_n(z)$ can then be described by
$$
T(x)\,\mapsto\,
\Ra_{\ts01}(\ts x\com z_1)
\ts\ldots\ts
\Ra_{\ts0n}(\ts x\com z_n)
\,\cdot\ts(\ts1\ot A_n)
\quad\qed
$$ 
\end{proof}

\noindent
Let us consider the one-dimensional
$\YN\ts$-modules $V_N(z)$ and $V_N(z)^\ast\ts$.
We will use the following observation, cf.\ \cite[Subsection 2.10]{MNO}.

\begin{Proposition}
In the $\YN\ts$-modules\/ $V_N(z)$ and $V_N(z)^\ast$, respectively
\begin{equation}\label{rat}
T_{ab}(x)\,\mapsto\,\ts\de_{ab}\,\cdot\,\frac{x-z-1}{x-z}
\,\quad\text{and}\quad\,\/
T_{ab}(x)\,\mapsto\,\ts\de_{ab}\,\cdot\,\frac{x-z}{x-z-1}\ .
\end{equation}
\end{Proposition}

\begin{proof}
It suffices to describe the action of $\YN$ only in the module $V_N(z)\ts$.
The decription of the action of $\YN$ in $V_N(z)^\ast$ will then follow
by the definition of the antipode $\operatorname{S}\ts$.
Using Proposition 2.4, the left hand side of the equality
(\ref{ext}) can be rewritten (\ts for any $n$) as
\begin{equation}\label{lhs}
\left(1-\frac{P_{\ts01}+\ldots+P_{\ts0n}}{x-z}\right)\,\cdot\ts(\ts1\ot A_n)\,\ts.
\end{equation}
If $n=N$, the sum $P_{\ts01}+\ldots+P_{\ts0n}$
can be replaced by $1$ without affecting the value of (\ref{lhs})\ts;
see for instance the proof of \cite[Theorem 2.3]{NT1}
\qed
\end{proof}

\noindent
Hence the one-dimensional $\YN\ts$-module, described in the display 
(\ref{rat}) first, is rational\ts; see Subsection 1.6. 
By applying Lemma 4.4 when $n=N$, the second $\YN\ts$-module in
(\ref{rat}) is also rational.
By using the definition (\ref{1.33}) of the
comultiplication on $\YN$, we now obtain the following result.

\begin{Corollary}
Let\/ $g(x)\in\CC(x)$ be a rational function such that\/ $g(\infty)=1$, with
all zeros and poles contained in $z+\ZZ$ for some common parameter $z\in\CC$.
Then the one-dimensional\/ $\YN\ts$-module defined by {\rm(\ref{onedimrep})},
is rational.
\end{Corollary}

\noindent\textbf{4.5.}
In this subsection, we shall derive Proposition 1.6 from Theorem 1.5.
The (\ts independent) proof of Theorem 1.5 will be given in the next
section. However, see the concluding remarks in Subsection 5.4.

Take the partition $\eta$ introduced in Subsection 1.6.
Using the definition (\ref{1.62}), we obtain the equality  
of rational functions
\begin{equation}\label{gid}
g_\eta(x+r)\,=\,
\frac{g_\mu(x)}
{g_{\ts\mut\ts}(\ts-\ts x\ns-\!M)}
\,\cdot\,
\frac{(\ts x+r)\ts(\ts x+M)}{x\ts(\ts x+M+r)}\ .
\end{equation}
Take the vector space $V_\xi^{\,\eta}$ defined as in
Subsection 1.2. Regard this vector space as $\YN\ts$-module
according to Subsection 1.5. The pullback of this module
relative to the automorphism $\tau_{-r}$ is denoted by  
$V_\xi^{\,\eta}(-r)\ts$.

By Theorem 1.5, the $\YN\ts$-module $\VOm$ 
is equivalent to the $\YN\ts$-module $\Vlm$, 
which was defined in Subsection 1.5. Using this definition,
the latter $\YN\ts$-module is equivalent to the pullback of the module
$V_\xi^{\,\eta}(-r)\ts$ relative to
the automorphism (\ref{1.61}) of $\YN\ts$, where
$$
g(x)\,=\,
\frac{g_\mu(x)\,\ts g_{\ts\lat}(\ts-\ts x\ns-\!M)}
{g_{\ts\mut\ts}(\ts-\ts x\ns-\!M)\,g_\eta(x+r)}
\cdot
\frac{(\ts x+r)}x
\,=\,g_{\ts\lat}(\ts-\ts x\ns-\!M)
\cdot
\frac{\,x+M+r}{x+M}\ .
$$
Here we used the definitions (\ref{1.52})\com(\ref{1.69})
and the equality (\ref{gid}).
Applying Theorem 1.5 to the $\YN\ts$-modules
$V_\xi^{\,\eta}$ and $V_\Gamma$ instead of $\Vlm$ and $\VOm$ respectively,
we now obtain Proposition 1.6.


\section{\hskip-3.5pt. Proof of Theorem 1.5}


\noindent\textbf{5.1.}
In this section, our basic vector space will be $\CC^L$ like
it was in Section~3, and unlike it was in the previous section.
Thus instead of (\ref{aux}),
we will work with the tensor product of vector spaces
\begin{equation}\label{aul}
\CC^L\ns\ot\Wll\ts=\ts
\CC^L\ns\ot(\CC^L)^\ast\ot\ldots\ot(\CC^L)^\ast\ot\CC^L\ns\ot\ldots\ot\CC^L
\end{equation}
where the first tensor factor $\CC^L$ will be auxiliary,
and shall not be regarded as a representation of $GL_L\ts$.
The $1\ns+\ts\lt+l\ts$ tensor factors in (\ref{aul}) will be labeled,
respectively from the left to right, by the indices (\ref{inda}).
For all possible indices $i\com j$ the operators $P_{ij}\com Q_{ij}$
and the functions 
$R_{ij}(x\com y)\ts$,
${}\,\Rb_{ij}(x\com y)\com\,\Ra_{0i}(x\com y)$
will be defined as their counterparts in Subsection~4.1, 
but $N$ will be replaced by $L\ts$.
In particular, these functions will take values~in
\begin{equation}\label{twoall}
\End(\CC^L)\ot\End(\Wll)=
\End(\ts\CC^L\ns\ot\Wll)\,.
\end{equation}
Instead of (\ref{instead}), we will use the second equality in
the display (\ref{3.555}). We will set
$L=N+M$, and fix a decomposition $\CC^{L}=\CC^N\ns\op\ts\CC^M$.

We begin the proof of Theorem 1.5 by choosing any standard tableau
$\La$ of shape $\la\ts$, such that the tableau $\Om$ is obtained from $\La$
by removing the boxes with numbers $1\lcd m\ts$. Here $\la$
is a partition of $l\ts$, and $m=l-n\ts$. Similarly, choose
any standard tableau $\Lat$ of shape $\lat\,$, such that the tableau $\Omt$
is obtained from $\Lat$ by removing the boxes with numbers $1\lcd\mt\ts$.
Again, here $\lat$ is a partition of $\lt\ts$, and $\mt=\lt-\nt\ts$.
The standard tableaux of non-skew shapes $\mu$ and $\mut\ts$, 
obtained by removing the boxes
with the numbers $m+1\lcd l$ from $\La\ts$
and the boxes with the numbers $\mt+1\lcd\lt$ from $\Lat\ts$,
will be denoted respectively by $\Up$ and $\Upt\ts$.

In Section 3, we denoted by $F_\La$ the linear operator
in the tensor product $W_l$ of $l$ copies of the vector space $\CC^L\ts$,
corresponding to the matrix element $f_\La\in\CC S_l\ts$.
Proposition 2.4 implies the equality of rational functions in~$x$
$$
R_{\ts0,\lt+1}(\ts x\com c_1(\La))
\,\ldots\,
R_{\ts0,\lt+l}(\ts x\com c_{\ts l}(\La))
\ts\cdot\ts(\ts1\ot1\ot F_\La)
$$
\vglue-12pt
\begin{equation}\label{4.61}
=\,
\biggl(1\ts-\!\sum_{\lt<j\le\lt+l}\frac{P_{\ts0j}}x\ts\biggr)
\cdot\ts(\ts1\ot1\ot F_\La)
\end{equation}
with values in the algebra (\ref{twoall}).
By applying Proposition 2.4 to the standard tableau $\Lat$ instead of
$\La\ts$, and by changing $x$ to $-x\ts$, we obtain the equality
$$
\Ra_{\ts01}(\ts x\com-c_{\ts\lt\ts}(\Lat))
\,\ldots\,
\Ra_{\ts0\ts\lt\ts}(\ts x\com -c_1(\Lat))
\ts\cdot\ts(\ts1\ot(\ts P_{\ts\lt}\,F_{\Lat}\,P_{\ts\lt}\ts)\ot1)
$$
\vglue-12pt
\begin{equation}\label{4.62}
=\ g_{\ts\lat}(-x)\ts\cdot
\biggl(1\ts+\!\sum_{1\le i\le\lt}\frac{Q_{\ts0i}}x\ts\biggr)
\cdot\ts(\ts1\ot(\ts P_{\ts\lt}\,F_{\Lat}\,P_{\ts\lt}\ts)\ot1)
\end{equation}
where $g_{\ts\lat}(-x)$ is defined via (\ref{1.62}). 
To obtain (\ref{4.62}), we have also used (\ref{fact}).

Observe that for any indices $i$ and $j$ such that 
$0<i\le\lt<j\ts$, we~have
$$
Q_{0i}\ts P_{0j}\ts\FLa\,=\,P_{0j}\,Q_{ij}\ts\FLa\,=\,0
$$
by Proposition 3.4. So the equalities
(\ref{4.61})\com(\ref{4.62}) and Lemma 1.3 imply~that
$$
\Ra_{\ts01}(\ts x\com-c_{\ts\lt\ts}(\Lat))
\ts\ldots\ts
\Ra_{\ts0\ts\lt\ts}(\ts x\com -c_1(\Lat))\,\ts
R_{\ts0,\lt+1}(\ts x\com c_1(\La))
\,\ldots\,
R_{\ts0,\lt+l}(\ts x\com c_{\ts l}(\La))
$$
\vglue-12pt
\begin{equation}\label{4.6}
\times\ \FLa\ =\ g_{\ts\lat}(-x)\ts\cdot
\biggl(1\ts+\!\sum_{1\le i\le\lt}\frac{Q_{\ts0i}}x
\ -\!\sum_{\lt<j\le\lt+l}\frac{P_{\ts0j}}x
\ts\biggr)
\cdot\ts\FLa\,\ts.
\end{equation}


\smallskip\noindent\textbf{5.2.}
The equality (\ref{4.6}) is the starting point for our proof of Theorem 1.5. 
Firstly, consider the case when $M=0$ and
$\mu=\mut=(0\com0\ts,\ts\ldots\ts)\ts$. In this case $N=L$, $n=l$, $\nt=\lt$
and the tableaux $\Om=\La$, $\Omt=\Lat$ have non-skew shapes.
In this case, the left hand side of the equality (\ref{4.6})
describes the action of the Yangian $\YL$ on the submodule $\VOm=\VLa$
of the tensor product of the dual evaluation modules with
the parameters $-\ts c_{\ts\lt}(\Lat)\lcd-\ts c_1(\Lat)$
and the evaluation modules with
the parameters $c_1(\La)\lcd c_{\ts l}(\La)\ts$.
The tensor factors have to be taken in the order as specified here,
see Proposition 4.1.

The right hand side of (\ref{4.6})
describes the action of $\YL$ in the module $\Vlm=\Vla$ as
defined in Subsection~1.5\ts; see Lemma 4.1.
Indeed, here we have $g_\mu(x)=g_{\ts\mut}(x)=1$ and $\pi_{NM}=\pi_N$. 
The image of the operator
$\FOm=\FLa\,$, as a $\glL\ts$-submodule in $\Wll\,$,
is equivalent to the $\glL\ts$-module $\Vla\ts$. 
Thus the equality (\ref{4.6})
implies Theorem 1.5 in the case when $M=0$.

Now suppose that $M\ge1$. 
Consider again the tensor product of the dual evaluation modules over $\YL$
with the parameters $-\ts c_{\ts\lt}(\Lat)\lcd-\ts c_1(\Lat)$ and 
the evaluation modules with the parameters $c_1(\La)\lcd c_{\ts l}(\La)\ts$.
Then consider the $\YL\ts$-module,
obtained by pulling this tensor product
back through the automorphism $\om_L$ of the algebra
$\YL\ts$; see (\ref{1.51}). The action of
$\YL$ in the latter module is described by the assignment 
$$
\sum_{a,b=1}^L\,E_{ab}\ot T_{ab}(x)\,\ts\mapsto\,
R_{\ts0,\lt+l}(\ts-x\com c_{\ts l}(\La))^{\ts-1}
\ldots\,
R_{\ts0,\lt+1}(\ts-x\com c_1(\La))^{\ts-1}
$$
\vglue-12pt
$$
\times\ \,
\Ra_{\ts0\ts\lt\ts}(\ts-x\com -c_1(\Lat))^{\ts-1}
\ts\ldots\,
\Ra_{\ts01}(\ts-x\com-c_{\ts\lt\ts}(\Lat))^{\ts-1}
$$
\vglue-6pt
$$
=\ \,
g_\la(-x)\ g_{\ts\lat}(x)^{-1}\ts g_{\ts\lat}(x\ns-\ns L)^{-1}
$$
\vglue-8pt
$$
\times\ \,
R_{\ts0,\lt+l}(\ts x\com-c_{\ts l}(\La))
\,\ldots\,
R_{\ts0,\lt+1}(\ts x\com-c_1(\La))
$$
\vglue-16pt
\begin{equation}\label{wmod}
\times\ \,
\Ra_{\ts0\ts\lt\ts}(\ts x\com c_1(\Lat)\ns+\ns L)
\,\ldots\,
\Ra_{\ts01}(\ts x\com c_{\ts\lt\ts}(\Lat)\ns+\ns L)\,;
\end{equation}

\medskip\noindent
see (\ref{fact}) and (\ref{3.55}),(\ref{3.6}),(\ref{3.75}).
Denote by $W$ the restriction of this $\YL\ts$-module
to the subalgebra $\YN\subset\YL\ts$;
we use the natural embedding $\ph_M:\YN\to\YL\ts$.
That is, $\ph_M:\ts T_{ab}(x)\ts\mapsto\ts T_{ab}(x)$ 
for $1\le a\com b\le N\ts$ by definition.
Further, denote by $\Wp$ the $\YN\ts$-module
obtained by pulling the $\YN\ts$-module $W$ back through
the automorphism $\om_N$ of $\YN\ts$.
Note that the vector space of the $\YN\ts$-modules
$W$ and $\Wp$ is $\Wll\,$. 

The subspace $\Zmm\subset\Wll$ is defined by (\ref{3.41}).
For every $i>\lt+m$ put
\begin{equation}\label{Pw}
R_{\ts0i}(x\com y)^{\,\wedge}\,=\,\ts
1-\frac{\,P_{\,0i}^{\ts\wedge}\,}{x-y}
\end{equation}
where $P_{\,0i}^{\ts\wedge}$ stands for the restriction
of the permutational operator $P_{0i}$ on the vector space
$\CC^L\ot\Wll$ to the subspace $\CC^N\ot\Zmm\subset\CC^L\ot\Wll\,$.
If $0<i\le\nt\/$, let $Q_{\ts0i}^{\,\wedge}$ be the linear operator
on $\CC^N\ot\Zmm\,$, acting as (\ref{1.45}) on the $i\ts$th tensor
factor $(\CC^N)^\ast$ and the $0\ts$th
tensor factor $\CC^N\ns$ of $\ts\CC^N\ot\Zmm$, and acting as the identity
in other tensor factors of $\CC^N\ot\Zmm\,$. Then put
$$
\Ra_{\ts0i}(x,y)^{\,\wedge}\,=\,\ts
\biggl(1+\frac{Q_{\ts0i}^{\,\wedge}}{x-y}\ts\biggr)
\cdot
\biggl(1-\frac1{(x-y)^{\ts2}}\biggr)^{\!-1}\ts.
$$

One can define an action of $\YN$ in the vector space $\Zmm$~by
$$
\sum_{a,b=1}^N\,E_{ab}\ot T_{ab}(x)\,\ts\mapsto\,
g_\la(-x)\ts\,g_{\ts\lat}(x)^{-1}\ts
g_{\ts\lat}(\ts x-\ns L)^{-1}\ts
g_{\ts\mut\ts}(\ts x-\ns L)
$$
\vglue-12pt
$$
\times\ \,
R_{\ts0,\lt+l}(\ts x\com-c_n(\Om))^{\,\wedge}
\,\ldots\,
R_{\ts0,\lt+m+1}(\ts x\com-c_1(\Om))^{\,\wedge}
\hskip-6pt
$$
\vglue-14pt
\begin{equation}\label{zmod}
\times\ \,
\Ra_{\ts0\ts\nt\ts}(\ts x\com c_1(\Omt)\ns+\ns L)^{\,\wedge}
\,\ldots\,
\Ra_{\ts01}(\ts x\com c_{\ts\nt}(\Omt)\ns+\ns L)^{\,\wedge}\,\ts.
\end{equation}

\smallskip\noindent
This assertion will follow from the observation below; let us 
denote by $Z$ the $\YN\ts$-module.
Further, let us denote by $\Zp$ the $\YN\ts$-module
obtained by pulling the $\YN\ts$-module $Z$ back through
the automorphism $\om_N$ of $\YN\ts$. 
The action of $\YN$ in $\Zp$ can be then described by the assignment
$$
\sum_{a,b=1}^N\,E_{ab}\ot T_{ab}(x)\,\ts\mapsto\,
g_{\ts\lat}(-x)\ts\,
g_{\ts\lat}(\ts-\ts x-\ns M)^{-1}\ts
g_\mu(x)^{-1}\ts
g_{\ts\mut\ts}(\ts-\ts x-\ns M)
$$
\vglue-12pt
$$
\times\ \,
\Ra_{\ts01}(\ts x\com-\ts c_{\ts\nt}(\Omt)\ns-\ns M)^{\,\wedge}
\,\ldots\,
\Ra_{\ts0\ts\nt\ts}(\ts x\com-\ts c_1(\Omt)\ns-\ns M)^{\,\wedge}
\hskip-24pt
$$
\vglue-8pt
$$
\times\ \,
R_{\ts0,\lt+m+1}(\ts x\com c_1(\Om))^{\,\wedge}
\,\ldots\,
R_{\ts0,\lt+l}(\ts x\com c_n(\Om))^{\,\wedge}\,\,.
$$

\smallskip\noindent
Observe that the $\YN\ts$-module $\Zp$
can also be obtained by pulling back the tensor product 
of $\YN$-modules
$$
V(\ts y_{\ts\nt})^\ast
\ot\ts\ldots\ot 
V(\ts y_1)^\ast
\ot\Wmmo\ot
V(\ts z_1)\ot\ts\ldots\ot V(z_n)
$$
through the automorphism (\ref{1.61}) of the algebra $\YN\ts$, where 
\begin{equation}\label{gx}
g(x)\,=\,
g_{\ts\lat}(-x)\,\ts
g_{\ts\lat}(\ts-\ts x-\ns M)^{-1}\ts
g_\mu(x)^{-1}\ts
g_{\ts\mut\ts}(\ts-\ts x-\ns M)\,.
\end{equation}
Here we keep to the settings (\ref{setz}),(\ref{sety}).
The Hopf algebra $\YN$ acts in the tensor factor $\Wmmo$ trivially,
that is via the counit homomorphism $\,\varepsilon\ts$.

The underlying vector space of the $\YN\ts$-modules
$Z$ and $\Zp$ is $\Zmm\,$, see (\ref{3.41}).
In Subsection 3.5, we introduced the projector
\begin{equation}\label{4.65}
H_{m\mt}:\,\Wll\,\to\,\Zmm\,.
\end{equation}

\begin{Proposition}
The projector\/ {\rm (\ref{4.65})} is\/
a\/ $\YN\ts$-intertwiner\/ $W\ns\to Z$.
\end{Proposition}

\begin{proof}
This proposition follows
by comparing the assignment (\ref{wmod})
with the definition (\ref{zmod}) of the $\YN\ts$-module $Z\ts$.
Indeed, let us expand the product 
at the right hand side of the equality (\ref{wmod}) as the sum
\begin{equation}\label{Sab}
\sum_{a,b=1}^L\,E_{ab}\ot S_{ab}(x)
\end{equation}
for certain functions
$S_{ab}(x)$ taking values in $\End(\Wll)\ts$.
It suffices to show that the sum
$$
\sum_{a,b=1}^N\,E_{ab}\ot(\ts H_{m\mt}\,S_{ab}(x))
$$
is equal to the product at the right hand side of the assignment
(\ref{zmod}), multiplied by $1\ot H_{m\mt}$ on the right.

To do this, let us apply to the values of the functions $S_{ab}(x)$
in the sum (\ref{Sab})
the operator conjugation with respect to each of the $\lt$ tensor
factors $(\CC^L)^\ast$ of $\Wll\,$, cf.\ the beginning of
our proof of Proposition 3.5. The resulting sum is the expansion
of the product of the rational functions in $x\ts$,
$$
g_\la(-x)\ts\,g_{\ts\lat}(x)^{-1}\,\cdot\ts
R_{\ts0,\lt+l}(\ts x\com-c_{\ts l}(\La))
\,\ldots\,
R_{\ts0,\lt+1}(\ts x\com-c_1(\La))
$$
\vglue-16pt
\begin{equation}\label{A}
\times\ \,
R_{\ts0\ts\lt\ts}(\ts-x\com-c_1(\Lat)\ns-\ns L)
\,\ldots\,
R_{\ts01}(\ts-x\com-c_{\ts\lt\ts}(\Lat)\ns-\ns L)
\end{equation}
taking values in the algebra of operators on
the tensor product 
of $1\ns+\lt+l$ copies
of the vector space $\CC^L$. Denote by $W_{1+\lt+l}$ this tensor product.
Each factor in (\ref{A}) is defined by
(\ref{3.45}) where $P_{0j}$ is the linear operator on $W_{1+\lt+l}$
exchanging the $0\ts$th and $j\ts$th copies of $\CC^L$.

The right hand side of (\ref{zmod}) is a rational function of $x$,
taking values in the operator algebra $\End(\CC^N\ot\Zmm\ts)\ts$.
By applying to these values the operator conjugation
with respect to each of the $\nt$ tensor
factors $(\CC^N)^\ast$ of $\CC^N\ot\Zmm\ts$, we obtain the function
$$
g_\la(-x)\ts\,g_{\ts\lat}(x)\,\cdot\ts
R_{\ts0,\lt+l}(\ts x\com-c_n(\Om))^{\,\wedge}
\,\ldots\,
R_{\ts0,\lt+m+1}(\ts x\com-c_1(\Om))^{\,\wedge}
$$
\vglue-12pt
\begin{equation}\label{B}
\times\ \,
R_{\ts0\ts\nt\ts}(\ts-x\com-c_1(\Omt)\ns-\ns L)^{\,\wedge}
\,\ldots\,
R_{\ts01}(\ts-x\com-c_{\ts\nt}(\Omt)\ns-\ns L)^{\,\wedge}
\end{equation}
taking values in the algebra of operators on
$\CC^N\ot W_{\nt}\ot\Wmmo\ot W_n\,$.
Each factor in (\ref{B}) is defined by (\ref{Pw})
where $P_{\,0i}^{\ts\wedge}$ stands for the restriction
of~the permutational operator $P_{0i}$ on the vector space 
$W_{1+\lt+l}$ to the subspace
$$
\CC^N\ot W_{\nt}\ot\Wmmo\ot W_n\,\subset\,W_{1+\lt+l}\ts\,.
$$
Proposition 5.2 follows by comparing
the products (\ref{A}) and (\ref{B}),
and by using Corollary 2.4
\qed
\end{proof}

\begin{Corollary}
The projector\/ {\rm (\ref{4.65})} is
a\/ $\YN\ts$-intertwiner\/ $\Wp\!\to\ts\Zp$.
\end{Corollary}


\noindent\textbf{5.3.}
Let us continue our proof of Theorem 1.5.
Consider the image of the subspace $\Zmm$
under the operator $\FLa$ in $\Wll\,$. Note that
this image
is contained in the subspace $\VLa\subset\Wll\,$.
Also consider the $\YN\ts$-module $W$ defined in Subsection 5.2\ts;
the vector space of this module is $\Wll\,$.

\begin{Proposition}
The image of the subspace\/ $\Zmm\subset\Wll$ under the operator $\FLa$
is a\/ $\YN$-submodule of\/ $W\ts$.
\end{Proposition}

\begin{proof}
The action of the coefficients of the series 
$T_{ab}(x)$ with $1\le a\com b\le N$
in the $\YN\ts$-module $W$ is described by the assignment (\ref{wmod}). 
Here we use the natural embedding $\ph_M:\YN\to\YL\ts$. 
Consider the product of the rational functions of $x$ at the
right hand side of (\ref{wmod})\ts; these functions take values in the
algebra (\ref{twoall}). We have an equality
$$
R_{\ts0,\lt+l}(\ts x\com-c_{\ts l}(\La))
\,\ldots\,
R_{\ts0,\lt+1}(\ts x\com-c_1(\La))
$$
\vglue-12pt
$$
\times\ \,
\Ra_{\ts0\ts\lt\ts}(\ts x\com c_1(\Lat)\ns+\ns L)
\,\ldots\,
\Ra_{\ts01}(\ts x\com c_{\ts\lt\ts}(\Lat)\ns+\ns L)
\,\cdot\,
(\,1\ot\FLa\ts)
\ =\ 
(\,1\ot\FLa\ts)
$$
\vglue-8pt
$$
\times\ \,
\Ra_{\ts01}(\ts x\com c_{\ts\lt\ts}(\Lat)\ns+\ns L)
\,\ldots\,
\Ra_{\ts0\ts\lt\ts}(\ts x\com c_1(\Lat)\ns+\ns L)
$$
\vglue-12pt
\begin{equation}\label{C}
\ts\ \times\ \,
R_{\ts0,\lt+1}(\ts x\com-c_1(\La))
\,\ldots\,
R_{\ts0,\lt+l}(\ts x\com-c_{\ts l}(\La))\,;
\end{equation}

\smallskip\noindent
see our proof of Proposition 1.4, which was given in Subsection 4.2.
To get the equality (\ref{C}), we also used the equality in (\ref{wmod}).
Let us expand the product of $\ts\lt+l\ts$ factors in the last two lines of
the display (\ref{C}), as the sum
$$
\sum_{a,b=1}^L\,E_{ab}\ot Y_{ab}(x)
$$
for certain functions $Y_{ab}(x)$ taking values in $\End(\Wll\ts)\ts$.
Then consider the restrictions of the operator values of the
functions $Y_{ab}(x)$ with $1\le a\com b\le N$
to the subspace $\Zmm\subset\Wll\,$. Arguing like in the proof
of Proposition 5.2, but using Lemma 2.4 itself instead of its
Corollary 2.4, we prove that
$$
\sum_{a,b=1}^N\,E_{ab}\ot (\ts Y_{ab}(x)\,|\,\Zmm)
\ = \ g_{\ts\mut\ts}(x\ns-\ns L)
$$
\vglue-8pt
$$
\times\ \,
\Ra_{\ts01}(\ts x\com c_{\ts\nt\ts}(\Omt)\ns+\ns L)^{\,\wedge}
\,\ldots\,
\Ra_{\ts0\ts\nt\ts}(\ts x\com c_1(\Omt)\ns+\ns L)^{\,\wedge}
$$
\vglue-8pt
$$
\ts\ \times\ \,
R_{\ts0,\lt+m+1}(\ts x\com-c_1(\Om))^{\,\wedge}
\,\ldots\,
R_{\ts0,\lt+l}(\ts x\com-c_n(\La))^{\,\wedge}\,\ts.
$$

\smallskip\noindent
In particular, the operator values of the functions $Y_{ab}(x)$
with $1\le a,b\le N$ preserve the subspace $\Zmm\,$.
Now Proposition~5.3 follows from (\ref{C})\qed
\end{proof}

\noindent
Since the $\YN\ts$-module $\Wp$ is obtained from $W$
by pulling back through an automorphism of $\YN\ts$,
Proposition 5.3 has the following corollary.

\begin{Corollary}
The image of the subspace\/ $\Zmm\subset\Wll$ under the operator $\FLa$
is a\/ $\YN$-submodule of $\Wp\ts$.
\end{Corollary}


\noindent\textbf{5.4.}
In this subsection we complete the proof of Theorem 1.5.
Let $V$ be the image of the subspace $\Zmm\subset\Wll$ 
under the linear operator
$$
H_{m\mt}\ts\FLa:\,\Wll\,\to\,\Zmm\,.
$$ 
By identifying $\Zmm$ with the tensor product 
$\Wmmo\ot W_{\ts\nt}^{\ts\ast}\ot W_n$ as we did in Subsection 3.5,
that is by 
exchanging the tensor factors $\Wmmo$ and $W_{\ts\nt}^{\ts\ast}$
of $\Zmm\ts$, we identify $V$ with a certain subspace of
$\Wmmo\ot W_{\ts\nt}^{\ts\ast}\ot W_n\,$. But due to the equality
(\ref{3.4444}), the latter subspace coincides with the subspace 
\begin{equation}\label{subspace}
V_{\ts\Up\Upt}\ot\VOm\,\subset\,\Wmmo\ot W_{\ts\nt}^{\ts\ast}\ot W_n\,.
\end{equation}

It follows from
Corollaries 5.2 and 5.3, that the subspace $V\subset\Zmm$
is a submodule in the $\YN\ts$-module $\Zp\ts$.
Let us now regard $V$ as $\YN\ts$-module, by using
the action of the algebra 
$\YN$ in $V$ inherited from $\Zp\ts$. Then
$V$ is a subquotient of the $\YN\ts$-module $\Wp$ by definition.

The observation on the $\YN\ts$-module $\Zp\ts$, made immediately
before Proposition 5.2, yields the following
description of the $\YN\ts$-module $V$. 
Take the $\YN\ts$-module $\VOm$ as defined in Subsection 1.4. 
Pull $\VOm$ back through the automorphism (\ref{1.61}) 
of $\YN\ts$, where $g(x)$ is given by (\ref{gx}).
Extend the resulting action of $\YN$ in 
$\VOm\ts$ to the subspace (\ref{subspace}) so that
$\YN$ acts in $V_{\ts\Up\Upt}$ trivially.
By exchanging the tensor factors $\Wmmo$ and $W_{\ts\nt}^{\ts\ast}$
of $\Wmmo\ot W_{\ts\nt}^{\ts\ast}\ot W_n$ in (\ref{subspace}),
we then obtain the $\YN\ts$-module $V$.

The subspace $V_{\ts\Up\Upt}\subset\Wmmo$ is equivalent to $V_{\ts\mu\mut}$
as a representation of the group $GL_M\ts$.
The subspace $\VLa\subset\Wll$ is equivalent to $V_{\la\lat}$ as a
representation of the group $GL_L\ts$.
Let us now regard $\VLa$ as $\YN\ts$-module
by pulling back through the homomorphism
$\pi_{NM}:\YN\to\UL\ts$, 
and then through the automorphism (\ref{1.61}) of $\YN$ where
$g(x)=g_{\ts\lat}(-x)\ts$.

\begin{Proposition}
$\YN$-module $V$ is a subquotient of\/ $\YN\ts$-module\/ $\VLa\ts$.
\end{Proposition}

\begin{proof}
By the definition (\ref{1.69}), we have
$$
\pi_{NM}=\,\pi_L\circ\ts\om_L\circ\ts\ph_M\circ\ts\om_N\ts.
$$
Consider $\VLa$ as a $\YL\ts$-module, using the definition from
Subsection 1.4. That is, the $\YL\ts$-module $\VLa$ is a submodule
of the tensor product of the dual evaluation $\YN\ts-$modules with
the parameters $-\ts c_{\ts\lt}(\Lat)\lcd-\ts c_1(\Lat)$
and the evaluation $\YN\ts$-modules with
the parameters $c_1(\La)\lcd c_{\ts l}(\La)\ts$.
The equality (\ref{4.6}) implies, that
the same $\YL\ts$-module can be obtained from the representation
$\VLa$ of the group $GL_L\ts$, by pulling back through the homomorphism
$\pi_{L}:\YL\to\UL\ts$, 
and then through the automorphism (\ref{1.61}) of $\YN$ where
$g(x)=g_{\ts\lat}(-x)\ts$.

Hence the $\YN\ts$-module $\VLa$ as defined immediately before stating
Proposition 5.4, 
can also be obtained by pulling the action of $\YL$ in $\VLa$
back through the homomorphism
$$
\om_L\circ\ts\ph_M\circ\ts\om_N:\YN\to\YL\,.
$$
So $\VLa$ is a submodule in the $\YN\ts$-module $\Wp$.
But by definition, $V$ is a quotient of a certain $\YN\ts$-submodule of
$\Wp$. The latter submodule of $\Wp$ is contained in $\VLa$\qed
\end{proof}

\noindent 
Consider the restriction of the representation $\VLa$ of 
the group $GL_L$
to the subgroup $GL_M\subset GL_L\ts$. 
Realize the vector space (\ref{1.0}) as
\begin{equation}\label{4.10}
{\rm Hom}_{\,GL_M}(\ts V_{\ts\Up\Upt}\ts\com\VLa\ts)\,.
\end{equation}
Since the image of the homomorphism $\pi_{NM}$ is contained in
the subalgebra of $GL_M\ts$-invariants $\AMN\subset\UL\ts$, the action
of the algebra $\YN$ in $\VLa$ from Proposition 5.4,
induces an action of $\YN$ in (\ref{4.10}).
This action of $\YN$ in the vector space (\ref{4.10})
is irreducible, see \cite[Section 2]{MO}.

The operator (\ref{4.65}) is $GL_N\times GL_M\ts$-equivariant,
and the vector space $V_{\ts\Up\Upt}\ot\VOm$ of the $\YN\ts$-module $V$
comes with a natural action 
of the groups $GL_N$ and $GL_M\ts$. The action of $GL_M$ in $V$
commutes with the action of the algebra $\YN\ts$.
By Proposition 5.4, the $\YN\ts$-module 
\begin{equation}\label{4.11}
{\rm Hom}_{\,GL_M}(\ts V_{\ts\Up\Upt}\ts\com V\ts)
\end{equation}
is a subquotient of 
(\ref{4.10}). The $\YN\ts$-module (\ref{4.10})
is irreducible, it must be equal to the $\YN\ts$-module
(\ref{4.11}). Here we used Proposition~3.6.

The $\YN\ts$-module (\ref{4.11})  
can also be obtained by pulling the $\YN\ts$-module $\VOm$
as defined in Subsection 1.5, back through the automorphism
(\ref{1.61}), where $g(x)$ is given by (\ref{gx}).
The proof of Theorem 1.5 is complete. 

Note that (\ref{4.11}) is also a subquotient of (\ref{4.10})
as a representation of the group $GL_N\ts$. Thus we obtain
Proposition 1.3 together with Theorem 1.5.

Let us make a few concluding remarks. One can prove Proposition~1.6
independently of Theorem 1.5, cf.\ \cite[Theorem A.4]{KW}. Then
one can derive Theorem 1.5 from its particular case $\nt=0\ts$, 
considered in \cite[Section 4]{N2}, by using Proposition 1.6.
We chose the present proof of Theorem 1.5, because
it generalizes the proof for $\nt=0\ts$.
Moreover, the method of our proof of Theorem 1.5 extends from
the group $GL_N$ to other classical Lie groups, 
the orthogonal group $O_N$ and the symplectic group $Sp_{\ts N}\,$;
see \cite[Section 5]{N2}.


\begin{acknowledgement}\hskip-9.5pt.
I am grateful to Grigori Olshanski
for numerous conversations. This work has been
supported by the EC under the
grant ERB-FMRX-CT97-0100.
\end{acknowledgement}




\begin{thebibliography}{MNO}

\bibitem[B]{B}
{R.\,Brauer},
\textit{On algebras which are connected with the
semisimple continuous groups},
{Ann.\ Math.}
\textbf{38}
(1937),
857--872.

\bibitem[C1]{C1}
{I.\,Cherednik},
\textit{On special bases of irreducible finite-dimensional representations
of the degenerate affine Hecke algebra},
{Funct.\ Analysis Appl.}
\textbf{20}
(1986),
87--89.

\bibitem[C2]{C2}
{I.\,Cherednik},
\textit{A new interpretation of Gelfand-Zetlin bases},
{Duke Math.\ J.}
\textbf{54}
(1987),
563--577.

\bibitem[CP]{CP}
{V.\,Chari and A.\,Pressley},
\textit{Fundamental representations of Yangians and
singularities of $R$-matrices},
{J.\ Reine Angew.\ Math.}
\textbf{417}
(1991),
87--128.

\bibitem[D]{D}
{J.\,Dixmier},
\textit{Alg\`ebres Enveloppantes},
Gauthier-Villars, Paris, 1974.

\bibitem[K]{K}
{K.\,Koike},
\textit{On the decomposition of tensor products of the representations
of the classical groups: by means of universal characters}, 
{Adv.\ Math.}
\textbf{74}
(1989),
57--86.

\bibitem[KW]{KW}
{R.\,King and T.\,Welsh},
\textit{Construction of $GL(n)$-modules using composite tableaux},
{Linear and Multilinear Algebra}
\textbf{34}
(1993),
99--122.

\bibitem[M]{M}
{I.\,Macdonald},
\textit{Symmetric Functions and Hall Polynomials},
Clarendon Press, Oxford, 1995.

\bibitem[MO]{MO}
{A.\,Molev and G.\,Olshanski},
\textit{Centralizer construction for twisted Yangians},
{Selecta Math.}
\textbf{6}
(2000),
269--317.

\bibitem[MNO]{MNO}
{\hskip2pt\hskip-4pt A.\,Molev, M.\,Nazarov and G.\,Olshanski},
\textit{Yangians and classical Lie algebras},
{Russian Math.\ Surveys}
\textbf{51} 
(1996),
205--282.

\bibitem[N1]{N1}
{M.\,Nazarov},
\textit{Yangians and Capelli identities},
{Amer.\ Math.\ Soc.\ Translations}
\textbf{181}
(1998),
139--163.

\bibitem[N2]{N2}
{M.\,Nazarov},
\textit{Representations of twisted
Yangians associated with skew Young diagrams},
{\tt math.RT/0207115\ts}.

\bibitem[NT1]{NT1}
{M.\,Nazarov and V.\,Tarasov},
\textit{On irreducibility of tensor products of Yangian modules},
{Internat.\ Math.\ Research Notices}
(1998),
125--150.

\bibitem[NT2]{NT2}
{M.\,Nazarov and V.\,Tarasov},
\textit{On irreducibility of tensor products of Yangian modules
associated with skew Young diagrams},
{Duke Math.\ J.}
\textbf{112}
(2002),
342--378.

\bibitem[O1]{O1}
{G.\,Olshanski},
\textit{Extension of the algebra $U(g)$ for infinite-dimensional classical
Lie algebras $g$, and the Yangians $Y(gl(m))$},
{Soviet Math.\ Dokl.}
\textbf{36}
(1988),
569--573.

\bibitem[O2]{O2}
{G.\,Olshanski},
\textit{Representations of infinite-dimensional classical groups,
limits of enveloping algebras, and Yangians},
{Adv. Soviet Math.}
\textbf{2} 
(1991), 
1--66.

\bibitem[S]{S}
{J.\,Stembridge},
\textit{Rational tableaux and the tensor algebra of\/ $gl_n\ts$},
{J.\ Comb.\ Theory}
\textbf{A46}
(1987),
79--120.

\bibitem[VK]{VK}
{A.\,Vershik and S.\,Kerov},
\textit{Characters and factor representations of the
infinite unitary group},
{Soviet Math.\ Dokl.}
\textbf{26}
(1982),
570--574.

\bibitem[W]{W}
{H.\,Weyl},
\textit{Classical Groups, their Invariants and Representations},
Princeton University Press, Princeton, 1946.

\bibitem[Y1]{Y1}
{A.\,Young},
\textit{On quantitative substitutional analysis I\ts} and \textit{II\/},
{Proc.\ London Math. Soc.}
\textbf{33}
(1901),
97--146  
and
\textbf{34}
(1902),
361--397.

\bibitem[Y2]{Y2}
{A.\ Young,}
\textit{On quantitative substitutional analysis VI\ts},
{Proc. London Math. Soc.}
\textbf{34}
(1932),
196--230.

\end{thebibliography}
\enddocument